\documentclass[a4paper,fleqn]{cas-sc}
\usepackage[figuresright]{rotating}
\usepackage{natbib,har2nat}
\usepackage{threeparttable}
\usepackage[titletoc]{appendix}
\usepackage{graphicx}
\usepackage{subfigure}
\usepackage{longtable}
\newcommand{\rnum}[1]{\uppercase\expandafter{\romannumeral #1\relax}}

\begin{document}
\let\WriteBookmarks\relax
\def\floatpagepagefraction{1}
\def\textpagefraction{.001}
\let\printorcid\relax    

\shorttitle{Collaborative electric vehicle routing with meet points}

\shortauthors{Zhou et~al.}

\title [mode = title]{Collaborative electric vehicle routing with meet points}

 \author[1]{Fangting Zhou}

 \ead{fangting@chalmers.se}

 \affiliation[1]{organization={Electrical Engineering, Chalmers University of Technology},
 city={Gothenburg},
 country={Sweden}}

\author[2]{Ala Arvidsson}
\ead{ala.arvidsson@chalmers.se}

\affiliation[2]{organization={Technology Management and Economics, Chalmers University of Technology},
city={Gothenburg},
country={Sweden}}

\author[3]{Jiaming Wu}
\ead{jiaming.wu@chalmers.se}

\affiliation[3]{organization={Architecture and Civil Engineering, Chalmers University of Technology},
    city={Gothenburg},
    country={Sweden}}

\author[1]{Balázs Kulcsár}
\ead{kulcsar@chalmers.se}

\begin{abstract}
In this paper, we develop a profit-sharing-based optimal routing mechanism to incentivize horizontal collaboration among urban goods distributors. This paper investigates a collaborative routing problem for urban logistics, in which the exchange of goods at meet points is optimally planned en route. We show that collaboration does not only reduce the total cost but also increases the profit of each company by sharing some customers and the related profit. Hence, we focus on solving a collaborative electric vehicle routing problem under constraints such as customer-specific time windows, opportunity charging, vehicle capacity, and meet-point synchronization. The proposed Collaborative Electric Vehicle Routing Problem with Meet Point (CoEVRPMP) is modeled as a mixed-integer nonlinear programming problem. We first present an exact method for optimal benchmarks via decomposition. To handle real-world problems, we suggest using a metaheuristic method: adaptive large neighborhood search with linear programming. The viability and scalability of the collaborative method are demonstrated via numerical case studies: (i) a real-world case of two grocery stores in the city of Gothenburg, Sweden, and (ii) a large-scale experiment with 500 customers. The results underline the importance of horizontal collaboration among delivery companies. Collaboration helps to reduce the environmental footprint (total energy consumed) and to increase the individual company's profit at the same time. 
\end{abstract}

\begin{keywords}
Urban logistics \sep Collaborative vehicle routing \sep Electric vehicle \sep Meet point \sep Profit sharing \sep MINLP
\end{keywords}

\maketitle

\section{Introduction}

Cities are continuously experiencing growing demand for freight transportation \citep{savelsbergh2016}. A 16\% annual growth rate in urban logistics is projected over the next five years from 2021, only connected to e-commerce \citep{Reuters2022}. Traffic congestion and greenhouse gas emissions are expected to increase by 21\% and 32\% until 2030, respectively \citep{World2020}. However, it has been demonstrated that delivery vehicles often operate below their capacity, delivering nothing more than "air" \citep{verlinde2012,chen2016}. The need for reliable and timely transportation solutions to balance the interests of society, businesses, and customers has never been more crucial.

In response to these challenges, horizontal collaboration through sharing economy business models, such as sharing logistics infrastructure and services with competitors, has emerged as a potential solution \citep{los2020,salama2022,DHL2022}, and it is gaining traction among practitioners and researchers \citep{pan2019,ferrell2020}. Such collaboration typically involves companies with shared interests and businesses. The majority of studies indicate that such collaboration can enhance the non-collaborative solution by approximately 20-30\% \cite{Gansterer2018}. However, it's important to note that existing studies often enforce collaboration from a holistic perspective, overlooking the individual benefits for each company. This may sacrifice a company in order to achieve larger total profits, discouraging horizontal collaboration in practice.

Moreover, to address sustainability development needs in the transportation sector, the adoption of electric vehicles (EVs) for goods distribution is gaining support as a viable solution \citep{malladi2022,yang2022}. In addition to advancements in vehicular technology and investments in charging infrastructure \citep{ghamami2020,mccabe2023}, the transition is impeded by route planning concerns regarding delivery range. The Electric Vehicle Routing Problem (EVRP), as seen in \citet{schneider2014}, \citet{keskin2016} and \citet{basso2019}, seeks to bridge the planning gap between limited range and effective urban distribution. Integrating electric vehicles into horizontal collaboration introduces new benefits but, at the same time, new challenges related to charging and route planning integration.

\begin{figure*} [htbp]
	\centering
		\includegraphics[scale=.7]{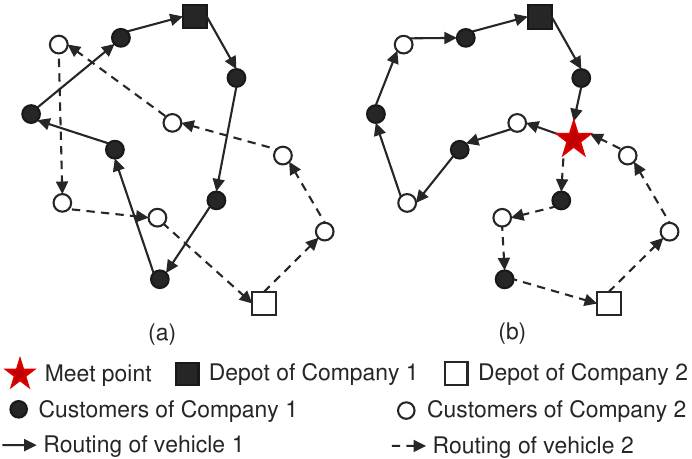}
    \caption{Horizontal collaboration example (a) non-collaboration, (b) collaboration}
	\label{FIG:1}
\end{figure*}

This paper introduces the concept of collaborative routing involving the exchange of goods en route at meet points. A visual representation of this idea is depicted in Fig. \ref{FIG:1}. The example involves two logistics companies serving their respective customers in the same area. Fig.\ref{FIG:1}(a) depicts the vehicle routes if company 1 (black) and company 2 (white) serve only their own customers. Fig.\ref{FIG:1}(b) demonstrates the collaboration scenario, where vehicles from both companies can exchange parcels at a "meet-point". In the case of collaboration, company A (B) serves not only its original customers but also the shared ones from the other company B (A). Each company offers three types of service: i) serving its own customers from the depot to the end; ii) serving its own customers from the depot to the meet point for exchange; iii) serving the other company's customers from the meet-point to the end. Due to the joint activities, a profit-sharing mechanism is introduced to split the profit from a shared customer, based on the profit ratio concept, further explained in Section \ref{Section3}.

This paper studies the collaborative electric vehicle routing problem with meet points (CoEVRPMP), explicitly considering individual companies' benefits. We explore a scenario where two logistics companies collaborate to plan vehicle routes to cross-serve a strategically selected set of customers. Individually serving these customers would be cost-prohibitive for either company. Instead, a unified global optimum solution is designed with the aim to increase the profitability of each individual company through collaboration and reduce the overall costs compared to non-collaborative solutions. We assume that the companies opt to transfer goods at several designated meet points and share customer addresses when collaborating (with standardized shipments). Various factors, such as customer-specific time windows, vehicle capacity, charging schedules, and meet-point synchronization, are taken into account. To address these challenges, we've developed a solution for CoEVRPMP, suitable for small to medium-sized real-world scenarios, using both exact and heuristic methods, and with the potential to scale for larger cases of up to 500 customers. The contributions of this paper can be summarized as follows:

\begin{itemize}
    \setlength{\itemsep}{0pt}
    \setlength{\parsep}{0pt}
    \setlength{\parskip}{0pt}
    \item The concept of meet points (transshipment points) is introduced to the collaborative routing problem, accompanied by a clear profit-sharing mechanism.
    \item The CoEVRPMP is formally defined and modeled as a mixed integer nonlinear programming problem.
    \item Practical constraints, including charging, customer time windows, vehicle capacity, and meet-point synchronization, are explicitly considered in an integrated framework. 
    \item An exact method and a metaheuristic algorithm are developed for theoretical analysis and practical implementation purposes, respectively.
\end{itemize}

The remainder of this paper is organized as follows. Section \ref{Section2} reviews the literature related to the CoEVRPMP. Section \ref{Section3} presents the problem description, formulates the mathematical programming model, and describes the proposed two solution methods: one exact and one metaheuristic. Section \ref{Section4} presents the experimental study and discusses the numerical results. Finally, Section \ref{Section5} concludes the paper with directions for future research.

\section{Literature review} \label{Section2}

The collaborative vehicle routing problem (CoVRP) is an operational planning challenge within horizontal collaboration \cite{Gansterer2018}. Most CoVRPs focus on either routing optimization \citep{sprenger2014,perez2015,montoya2016,quintero2016,stellingwerf2018,munoz2019,vahedi2022} or profit sharing \citep{berger2010,curiel2013}. However, only a limited number of studies have addressed both aspects \citep{krajewska2008,zibaei2016,wang2017}. A more comprehensive review of collaborative vehicle routing can be found at \citet{Gansterer2018}. Another critical consideration is the impact of electric vehicles on collaborative routing. This section provides a comprehensive review of these three aspects: routing optimization in CoVRPs, profit sharing in CoVRPs, and the integration of electric vehicles into CoVRPs.

\subsection*{Routing optimization in CoVRPs}

\medskip

\subsubsection*{Collaborative routing Paradigms}

Collaborative vehicle routing primarily falls into two categories: centralized planning and decentralized planning. Unlike decentralized planning, which entails limited or no information exchange, centralized planning involves information sharing. Centralized collaborative planning prioritizes optimizing the entire system over individual companies, while decentralized planning emphasizes more localized and independent decision-making.
Additionally, within the literature, there exists a distinction between two types of customer requests: 'reserved' and 'shared.' Reserved requests pertain to customers whom carriers must serve due to contractual obligations or other specific considerations, while shared requests encompass those customers whom carriers are open to serving collaboratively with others.

Centralized collaborative planning studies assess the potential benefits of collaborative versus non-collaborative settings. The potential benefits could be based on total costs \citep{lin2008}, total travel distance \citep{montoya2016,perez2015}, profits \citep{li2016,fernandez2016}, and emissions \citep{perez2015}. However, centralized collaborative planning focuses more on the whole system than the single company. Hence, one possible breakthrough is to incorporate individual profit gains into centralized collaborative planning.

There is limited research that has focused on centralized collaborative planning with profit gains and reserved customers. \citet{fernandez2016} propose a collaborative uncapacitated arc routing problem with profit gains, where the goal is to maximize the total profit of the coalition of carriers and take the lower bound on the individual profit of each carrier into account. The model considers side payments for those customers that are served by different carriers. Their work is based on the arc routing problem that sets arc as a customer, and the time windows of customers are ignored. Additionally, the side payments are hard to set. Thus, in this paper, reasonable side payments are set for each request.

In most of the centralized collaborative routing literature, depots could directly serve other companies' customers, where a strong assumption is that the collaborating companies share the same depot \citep{stellingwerf2018} or multiple depots \citep{perez2015,montoya2016,quintero2016,munoz2019}. Some studies consider exchanging goods between depots \citep{sprenger2014,vahedi2022}, which reduces the total cost but also brings additional travel costs to connect depots. Consequently, most of the centralized problems are formulated as the VRP or multi-depot VRP (MDVRP), but those problems have been widely studied. The main difference between the proposed non-collaborative and collaborative routing problems occurs only in the customer sets in these studies. \citet{juan2014} associate collaborative routing with backhaul, which is a simple collaborative method that merges two routes from different companies to reduce backhaul. In this way, the merged route visits customers after visiting their depot. Pickup and delivery (PD) requests are frequently added extensions here \citep{krajewska2008,wang2014,buijs2016,li2016}, where PD locations do not coincide with a depot. Then, requests are served and fulfilled before the vehicle returns to the depot. Thus, the depot is no longer needed to store goods, let alone to share depots or exchange goods. Regardless of whether depots are shared or connected via routes, the issue that needs to be addressed in this study is how to achieve better collaborative distribution by exchanging goods among companies.

Considering that carriers are often resistant to sharing all their customers' data with a central planner, some research has focused on decentralized planning, including request selection and request exchange \citep{Gansterer2018}. Regarding the request selection method, carriers need to decide which of their customers can be offered to the collaboration partners. This is essential because some companies may not be willing to share all of their customers. The exchange of goods for customers could be included in vehicle routes (lane exchanges) or via auction-based systems. However, the sharing preferences of collaborators limit significant profit increases. An interesting decentralized planning study by \citet{li2016} proposes a pickup and delivery problem (PDP) with time windows, profits, and reserved customers in carrier collaboration realized through combinatorial auction. This research focuses only on one carrier and includes two decisions: which customers to bid for (to serve) and how to build routes for maximizing "own profit". Like many of the decentralized planning studies, \citet{li2016} have a myopic focus: increase profit share for a single company. Therefore, a valid research question is raised on how to jointly ensure the companies’ profit while lowering the total cost of the whole system.

\subsubsection*{Transshipment in CoVRPs}

Only a few papers have studied the routing problems with transshipment, such as PDP with transshipment (PDPT) \citep{cortes2010}, Vehicle Routing Problem with Transhipment Facilities \citep{baldacci2017}, and Two-Echelon Vehicle
Routing Problem \citep{crainic2009}. \citet{mitrovic2006} assess the usefulness of transshipment and state that transshipment points prove highly beneficial in clustered instances. \citet{drexl2012} emphasized critical challenges in addressing synchronization aspects, including the PDPT and its related problem variations. Research shows that the benefit of allowing transshipment can be significant \citep{lyu2023}. The transshipment in the above studies is within a single company. Expanding the concept of transshipment among companies may enhance collaboration and yield further benefits. This is one of the objectives of this paper.

The closest study \citet{zhang2022} addressed the goods exchange issue by transferring goods at customer points or depots and studied a heterogeneous multi-depot collaborative vehicle routing problem. This work shows that transferring goods en route (from unloading vehicle to loading vehicle) can result in different gains in the system. However, several aspects can be added to this study to increase its practical applicability and relevance, which will be addressed in this paper. These aspects include i) exchanging between two vehicles instead of only from unloading vehicle to loading one; ii) time windows of customers and specified waiting time at transfer points; iii) profit sharing or minimal profit guarantee for the initiatives of collaboration.

\subsection*{Profit sharing in CoVRPs}

An important aspect of collaborative operations is how to share the potential extra profit among the collaborators. This calls for the solution of cost allocation problems \cite{engevall2004}. \citet{guajardo2016} review cost allocation solutions for collaborative transport services and summarize the most commonly used methods. This includes the commonly used Shapley value \citep{vanovermeire2014,kimms2016} and other proportional methods \citep{berger2010,ozener2013}. Note that these methods of sharing profit require knowing the total benefit first.

Only a few studies integrate routing planning with profit-sharing aspects in the design of collaborative vehicle routing problems. \citet{krajewska2008} combine routing and scheduling problems with cooperative game theory. It proposes two subproblems to be addressed and integrated. First, it hints at solving the routing problem (multi-depot PDP with time windows). Second, a profit-sharing mechanism involves the Shapley value to determine a fair allocation. However, profit sharing of this type may have potential legal risks, e.g., against antitrust or competition laws.

\subsection*{EVs Integration in CoVRPs}

The electric vehicle routing problem (EVRP) emerged from the traditional VRP by considering battery constraints, charging operations, and energy consumption. One of the earliest works on EVRP has been communicated in \citet{conrad2011}. It introduces the recharging vehicle routing problem, where vehicles with limited range are allowed to recharge at customers' locations. The recharging time is assumed to be fixed. \citet{schneider2014} study the electric vehicle routing problem with time windows and recharging stations (EVRPTW). The EVRPTW considers customer time windows and includes the possibility of opportunity recharging at stations, with the recharging time being dependent on the battery level. A comprehensive review of EVRP can be referred to \citet{kucukoglu2021}, where EVRP studies are classified according to four criteria: objective function types, energy consumption computations, considered constraints in the EVRP, and fleet types.

Very few studies incorporate the collaborative strategy in the EVRP. \citet{munoz2019} assess the implementation of an electric fleet of vehicles in urban goods distribution under a horizontal collaboration strategy between carriers. A multi-objective optimization is proposed in their study to explore the relationship between the delivery cost and the environmental impact. \citet{vahedi2022} study a collaborative capacitated electric vehicle routing problem, where a bi-objective function is considered to minimize i) the total tardiness costs and fixed costs of using EVs and ii) the total electrical energy consumption. They assume that there is a Third Party Logistics company to transship goods between depots and include this cost in the objective function. The above two studies model the collaborative scenario as MDVRP and multi-depot EVRP, respectively. The former is similar to most of the centralized planning studies, and the latter injects the electric vehicle characters into the MDVRP setup. Both of these studies overlooked the profit of the individual company and did not consider the time windows of customers and reserved customers of companies.

Unlike transferring products unilaterally from one vehicle to another, as seen in \citet{zhang2022}, this paper emphasizes bilateral exchanges between two vehicles, which may operationally be more tractable. Moreover, we propose an optimization-driven mechanism to exchange goods en route for collaboration, as opposed to depot-based transfers \citep{sprenger2014,perez2015,montoya2016,munoz2019,vahedi2022}. Regarding electric vehicles, only one study \citep{vahedi2022} has expanded the basic collaborative routing by adding charging possibilities, which is usually formulated as MDVRP that is not directly related to collaboration.

This paper investigates collaborative electric vehicle routing problems within a centralized planning framework, focusing on collaboration and electric vehicles (EVs). The collaboration involves the exchange of goods and profit-sharing, with a specific focus on partial EV charging. In contrast to existing literature, we consider scenarios where vehicles can exchange goods en route. Through profit-sharing, we aim to reconcile conflicts between system-wide optimization and individual benefits. Notably, our approach integrates route optimization and profit-sharing within a comprehensive structure, seamlessly incorporating profit-sharing into the optimization process. Consequently, our model allows for the simultaneous derivation of optimal routing and profit-sharing solutions. In doing so, this study addresses critical practical constraints, including charging challenges, time windows, vehicle capacity, and synchronization at meet points.

\section{Methodology} \label{Section3}

In this paper, carriers collaborate by exchanging goods at one of several designated 'meet points'. This interaction occurs because their delivery routes intersect, presenting significant decision-making challenges, including selecting the meet points, ensuring vehicle arrivals are synchronized at these points, and integrating them into route optimization. Additionally, we address issues like en-route charging and customer-specified time windows, which intricately link vehicle routes, energy consumption, and partial charging strategies. These considerations contribute to the complexity of the collaborative routing problem.

Without loss of generality, the following assumptions (boundary conditions \textbf{A2}-\textbf{A5}, model/method specific assumptions \textbf{A1}, \textbf{A6}-\textbf{A9}) are used along the paper. 
\begin{enumerate}
    \item[\textbf{A1}] Two companies are considered with one electric vehicle each, starting from and returning to the same depot. \footnote{The proposed model can also be applied to multiple vehicles with slight modifications, which can be found in Appendix A. For ease of communication, we focus on the two-vehicle case in the main text.}
    \item[\textbf{A2}] Each company has two known sets of customers: a set of reserved customers to be served only by the company itself (due to company policy, privacy, user agreements, etc.) and a set of customers to share for collaboration.
    \item[\textbf{A3}] Each company has certain expectations for the profits of collaboration. Thus, the profit threshold will be defined by each company separately (based on strategic purpose, long-term development, etc.), below which companies will \emph{refuse} to collaborate. Since the companies' expectations are different from case to case, we deem it irrelevant to this study. In this work, we simply define the threshold as the non-collaborative profit (maximum profit achieved by a company operated independently).
    \item[\textbf{A4}] There exists a mutually trusted consolidator. The collaboration is planned in a centralized manner, which means their information should be provided to the central planner, and both companies comply if agreed.
    \item[\textbf{A5}] Electric vehicles can be put on charge at customer locations and at meet-points, where partial charging is considered.
    \item[\textbf{A6}] Electric vehicles are fully charged when leaving the depot.
    \item[\textbf{A7}] Electric vehicle capacities are deterministic and known.
    \item[\textbf{A8}] Each customer is visited by only one company, but the full chain of service may involve another company if they exchange goods at meet points.
    \item[\textbf{A9}] The travel time, the delivery time window, and the travel distance among customers are known to be deterministic.
\end{enumerate} 

With the above assumptions, we study the CoEVRPMP with predefined profit thresholds, time windows, state of charge and charging constraints, vehicle capacity, and meet-point synchronization. In the CoEVRPMP, we optimize several vital decisions to minimize total collaborative operational costs. These decisions encompass meeting time and location, assignment of the shared customer, vehicle delivery sequence, charging locations, and the amount of energy to charge. This section provides an overview of the optimization model and the solution approaches.

\subsection{Model formulation} \label{Section3.1}

To help the reader understand the CoEVRPMP, we now provide a mixed-integer nonlinear programming (MINLP) formulation of the problem. The CoEVRPMP is modeled using a complete directed graph $G=(N, A)$, where $N=O \cup R \cup M$ represents the node set and $A$ is the edge set. Specifically, $O$ is the depot set, $R$ represents the customer set, and $M$ is the meet point set. The customer set $R$ comprises two subsets: reserved customers $R^r$ and shared customers $R^s$. Moreover, each company $k$ possesses a set of customers $R_k$, which can be further divided into reserved customers $R_k^r$ and shared customers $R_k^s$, with $k$ belonging to the company (vehicle) set $K=\left\{1,2\right\}$. 

The MINLP uses the following decision variables. Binary variables $x_{ij}^k$ take value 1 if vehicle $k$ delivers from node $i$ to node $j$. Binary variables $y_j^k$ take value 1 if customer $j$ is served by vehicle $k$. Binary variables $\varepsilon_m$ take value 1 if vehicles choose to meet at meet point $m$. Binary variables $z_i^k$ take value 1 if vehicle $k$ charges at node $i$. Variables ${\Phi}_k$ refer to the total profit of company $k$ (SEK). Variables $b_i^k$ and ${\delta}_i^k$ specify the remaining energy and the amount of battery charged for vehicle $k$ at node $i$ (Wh). Variables ${ST}_i^k$ are the time for serving goods and charging of vehicle $k$ at node $i$. Variables $s_i^k$ represent the service start time of vehicle $k$ at node $i$. Additionally, we introduce variables $T^k$ to denote the arrival time of vehicle $k$ at the end depot, which aligns with the service start time at node $i$ for vehicle $k$, where $i$ corresponds to the end depot and is within the set $O^-$.

For the convenience of communication, all the notations used for problem formulation are presented in Appendix B (Table \ref{tbnotation}). In addition, if a customer belongs to Company 1 but is also partially served by Company 2, we refer to Company 1 as the responsible company and Company 2 as the collaborative company, and vice versa. The MINLP formulation follows.

The profit is determined by subtracting the total delivery cost from customer revenue. Given constant customer revenue, lowering the overall delivery cost directly boosts profit. The objective is to minimize the total cost for all companies, leading to profit-maximizing:

\vspace{-\topsep}
\begin{equation} \label{Obj}
	\mathrm{min}\;\;\sum_{k\in K}\sum_{i\in N}\sum_{j\in N}c_d{D_{ij}}{x_{ij}^k} + \sum_{k\in K} c_t T^k,
\end{equation}

\noindent where the first term signifies the energy consumption cost associated with the distance $D_{ij}$, while the second term represents labor cost tied to arrival time. $c_d$ stands for unit energy consumption cost, and $c_t$ represents unit driver salary. Without loss of generality, the following equality and inequality constraints are defined. 

\vspace{2 mm}
\textbf{\textit{(\rnum{1}) Profit threshold constraints}}

In practice, collaboration can be highly motivated by a win-win situation, which in this context, means an increase in profit for both companies. Therefore, to make the results meaningful and practical, we introduce a profit-sharing threshold as a necessary condition of collaboration. Each company could determine its own threshold $P_k^{min}$, and the collaboration will only occur if the profit $\Phi_k$ exceeds the threshold $P_k^{min}$ for both companies, which can be formulated as:

\vspace{-\topsep}
\begin{equation} \label{NonPk}
	{\Phi}_k \geq P_k^{min}, \forall k \in K,
\end{equation}

\noindent where 

\vspace{-\topsep}
\begin{equation} \label{Profitk}
	{\Phi}_k = \sum_{j\in R_k}p_j{y_j^k} + \sum_{m\in M}\sum_{j\in R_k}p_j\alpha_j^{m}\varepsilon_m \left(1-y_j^k\right) + \sum_{m\in M}\sum_{j\in R\setminus R_k}p_j\left(1-\alpha_j^{m}\right)\varepsilon_m y_j^k - \sum_{i\in N}\sum_{j\in N}c_d{D_{ij}}{x_{ij}^k} - c_t T^k,
\end{equation}

\vspace{-\topsep}
\begin{equation} \label{Calpha}
\alpha_j^{m} = \frac{D_{{o_k}m}}{D_{{o_k}m}+D_{mj}}, \forall{j\in R_k},{k\in K},{m\in M}.
\end{equation}

The profit of a company $\Phi_k$ is naturally defined as the net income (income deducting cost) in Eq. \eqref{Profitk}, with the service fee of customer $j$ represented as $p_j$. The income comes from providing service to the three categories of customers, as corresponding to the first three terms in the equation, respectively. Specifically, the first term denotes income from customers entirely served by the responsible company; the second term represents income from customers partially served by the responsible company; and the third term accounts for the income from shared customers of the other company. Clearly, there is a need for a profit-sharing mechanism to split the income from shared customers.

The profit ratio $\alpha_j^m$ in Eq. \eqref{Calpha} serves as the core of our profit-sharing mechanism, which is a distance-based approach. With the defined ratio, we provide more insights into ratio penalized terms of income function in Eq. \eqref{Profitk}. If the two companies jointly serve customer $j$, the profits of the responsible company (serving from depot to meet point) and the collaborative company (serving from meet point to customer) are $p_j\alpha_j^{m}$ and $p_j\left(1-\alpha_j^{m}\right)$, respectively. While enabling the split of income, the profit-sharing mechanism introduces the complex interplay between the selection of meet points and shared customers (i.e., the multiplication of $\varepsilon_m$ and $y_j^k$ in Eq. \eqref{Profitk} ). This interplay makes the optimization model nonlinear and thus computationally intensive (see more details in Section \ref{SolutionApproach}). Last but not least, the remaining two terms in Eq. \eqref{Profitk} are the energy consumption and labor cost, respectively.

\vspace{2 mm}
\textbf{\textit{(\rnum{2}) Charging and capacity constraints}}

We now ensure that the delivery vehicles are running under practical capacity and favorable battery levels. In existing studies, it has been found that a high depth of discharge
exacerbates battery degradation \citep{schoch2018enhancing}. Thus, it is beneficial for companies to regulate the battery level of EVs and charge it en route. To this end, we constrain the EV battery energy within a lower and upper bound $[L, B]$, as follows, and enable charge:

\vspace{-\topsep}
\begin{equation} \label{BatteryCons}
	L \leq b_i^k \leq B, \forall{i\in N},{k\in K},
\end{equation}

\noindent while visiting customers, the battery state is updated by Eq. \eqref{BatteryLevel}, and opportunity charging is regulated in Eq. \eqref{chargeb}. Notably, energy consumption is directly linked to travel distance, with $\epsilon$ denoting the unit energy consumption per distance. 

\vspace{-\topsep}
\begin{equation} \label{BatteryLevel}
	b_j^k \leq b_i^k + {\delta}_i^k - \epsilon D_{ij} + B \left(1 - x_{ij}^k \right), \forall{i\in N},{j\in N},{k\in K},
\end{equation}

\vspace{-\topsep}
\begin{equation} \label{chargeb}
	{\delta}_i^k \leq \left(B - b_i^k\right)z_i^k, \forall{i\in N},{k\in K}.
\end{equation}

Eq. \eqref{capacityV} further ensures that the overall demands of the customers to be visited (where $q_j$ denotes the demand of customer $j$), encompassing both own and other customers' demands, do not exceed the capacity $Q_k$ of vehicle $k$.

\vspace{-\topsep}
\begin{equation} \label{capacityV}
	\sum_{i\in N}\sum_{j\in R} q_j x_{ij}^k \leq Q_k, \forall{k\in K}.
\end{equation}

\vspace{2 mm}
\textbf{\textit{(\rnum{3}) Time window constraints}}

Exchanging goods at the meet-point, the fundamental enabler of collaboration, entails space and time synchronization between the two vehicles in terms of their arrival time at the meet-point. In our study, a maximum waiting time window ${WT}_{max}$ is predetermined to ensure the vehicles can meet each other: 

\vspace{-\topsep}
\begin{equation} \label{MaxWT}
	\left|s_{m}^1 - s_{m}^2\right| \leq {WT_{max}}, \forall{m\in M}.
\end{equation}

In the time domain, the following constraints are further defined to ensure the vehicles deliver goods within the desired time windows of customers:
\vspace{-\topsep}
\begin{equation} \label{ServiceSeq}
	s_{m}^k - \Gamma \left(1-y_j^k\right) \leq s_j^k, \forall{j\in R^{\rm{s}}-R_k^{\rm{s}}},{k\in K},{m\in M},
\end{equation}

\vspace{-\topsep}
\begin{equation} \label{StartServiceTimeN}
	s_i^k + {ST}_i^k + {tt}_{ij}^k - \Gamma\left(1-x_{ij}^k\right) \leq s_j^k, \forall{i\in N},{j\in N},{k\in K},
\end{equation}

\vspace{-\topsep}
\begin{equation} \label{calculateST}
	{ST}_i^k = {st}_i + 60{\delta}_i^k/{r_i}, \forall{i\in N},{k\in K},
\end{equation}

\vspace{-\topsep}
\begin{equation} \label{Arrivaltime}
	T^k = s_i^k, \forall{i\in O^-},{k\in K},
\end{equation}

\noindent where Eq. \eqref{ServiceSeq} guarantees that the exchanged goods must be delivered after the meet-point, and arrival time and dwell time (including charging time and service time ${st}_i$) at each customer are updated and regulated in Eq. \eqref{StartServiceTimeN} and Eq. \eqref{calculateST}. The travel time ${tt}_{ij}^k$ and dwell time $ST_i^k$ are utilized to compute the arrival time. Charging time at node $i$ is computed based on the amount of battery charged $\delta_i^k$ and charging rate $r_i$. Eq. \eqref{Arrivaltime} ensures that the arrival time of vehicle $k$ equals the start service time at the end depot. Customer time windows $\left[e_j,l_j\right]$ are ensured by: 

\vspace{-\topsep}
\begin{equation} \label{TW}
	e_j \leq s_j^k \leq l_j, \forall{j\in R},{k\in K}.
\end{equation}

\vspace{2 mm}
\textbf{\textit{(\rnum{4}) Route constraints}}

The route constraints make sure that each and every customer will be served only once, and the two vehicles will meet one time at the same meet-point:

\vspace{-\topsep}
\begin{equation} \label{CusVisitOnce}
	\sum_{k\in K}\sum_{i\in N}{x_{ij}^k} = 1, \forall{j\in R},
\end{equation}

\vspace{-\topsep}
\begin{equation} \label{OneMP}
	\sum_{i\in N}\sum_{m\in M}{x_{im}^k} = 1, \forall{k\in K},
\end{equation}

\vspace{-\topsep}
\begin{equation} \label{SameMP}
	\sum_{i\in N}{x_{im}^1} - \sum_{i\in N}{x_{im}^2} = 0, \forall{m\in M},
\end{equation} 

\vspace{-\topsep}
\begin{equation} \label{depotSandE}
	\sum_{j\in R\cup M}{x_{{o_k^+}j}^k} = 1, \forall{k\in K},
\end{equation}

\vspace{-\topsep}
\begin{equation} \label{depotSandE2}
	\sum_{i\in R\cup M}{x_{i{o_k^-}}^{k}} = 1, \forall{k\in K},
\end{equation}

\vspace{-\topsep}
\begin{equation} \label{reservedCus}
	\sum_{j\in N}{x_{ij}^k} = 1, \forall{i\in R_k^{\rm{r}}},{k\in K},
\end{equation}

\noindent where Eq. \eqref{CusVisitOnce} guarantees that all customers will be visited exactly once, Eq. \eqref{OneMP} ensures that each vehicle visits only one meet-point, and Eq. \eqref{SameMP} guarantees that both vehicles will visit the same meet-point. Eq. \eqref{depotSandE} and Eq. \eqref{depotSandE2} ensure that vehicle $k$ must start from and return to the depot $o_k$. Eq. \eqref{reservedCus} guarantees that reserved customers will be served by the responsive company.

\vspace{2 mm}
\textbf{\textit{(\rnum{5}) Flow conservation constraints}}

\vspace{-\topsep}
\begin{equation} \label{mandx}
	\varepsilon_m = \sum_{i\in N}{x_{im}^k}, \forall{k\in K},{m\in M},
\end{equation}

\vspace{-\topsep}
\begin{equation} \label{yandx}
	y_j^k = \sum_{i\in N}{x_{ij}^k}, \forall{j\in R},{k\in K},
\end{equation}

\vspace{-\topsep}
\begin{equation} \label{conservation}
	\sum_{i\in N}{x_{ij}^k} - \sum_{i\in N}{x_{ji}^k} = 0, \forall{j\in R\cup M},{k\in K},
\end{equation}

\noindent where Eq. \eqref{mandx} ensures if meet point $m$ is chosen, then vehicles must visit $m$, Eq. \eqref{yandx} indicates whether request $j$ is served by vehicle $k$ through the link $i-j$, and the conservation of the arriving and the departing vehicle at each node is ensured by the Eq. \eqref{conservation}.

\vspace{2 mm}
\textbf{\textit{(\rnum{6}) Decision variables and their domains}}

\vspace{-\topsep}
\begin{equation} \label{decV1}
	x_{ij}^k,y_{j}^k,z_{i}^k\in \left\{ 0,1 \right\}, \forall{i\in N},{j\in N},{k\in K},
\end{equation}

\vspace{-\topsep}
\begin{equation} \label{decV2}
	s_{i}^k,b_{i}^k,{\delta}_{i}^k,{ST}_{i}^k \geq 0, \forall{i\in N},{k\in K},
\end{equation}

\vspace{-\topsep}
\begin{equation} \label{varepsilonm}
	\varepsilon_m \in \left\{ 0,1 \right\}, \forall{m\in M},
\end{equation}

\vspace{-\topsep}
\begin{equation} \label{Profit}
	{\Phi}_{k}\geq 0.
\end{equation}

Lastly, the decision variables are injected via the Eq. \eqref{decV1}, \eqref{decV2}, \eqref{varepsilonm}, and \eqref{Profit}. Even though the cost function is linear in the decision variables, the proposed CoEVRPMP is nonlinear due to Eq. \eqref{NonPk}, \eqref{Profitk}, \eqref{chargeb}, and \eqref{MaxWT}. Finally, a mixture of real-valued and integer-valued decision variables is used.

\subsection{Solution approach} \label{SolutionApproach}

The formulated MINLP is computationally intensive owing to the nonlinearity, integer variables, and a large set of decision variables and hard constraints. We first develop an exact approach that could theoretically achieve global optimality through linearization techniques. However, for large-scale problems (50 customers and up), the exact method may be computationally intractable with computer capacity nowadays (as shown later in Section \ref{Section4}). Therefore, we develop a metaheuristic method to facilitate real-world implementation.

\subsubsection{Exact algorithm via branching}\label{branchMILP}

Although the MINLP problem is generally undecidable \cite{jeroslow1973}, we notice that the nonlinearity of our problem is mainly due to Eq. \eqref{Profitk}. The structure paves the way to reduce the numerical complexity from undecidable to sequentially NP-hard. We present a technique to linearize the nonlinearity injected by the term $\varepsilon_m y_j^k$. Given that both decision variables are binary, one approach to linearize the term is by introducing an additional binary decision variable, denoted as $\Lambda = \varepsilon_m y_j^k$. This new decision variable, $\Lambda$, must adhere to the following conditions: 1) $\Lambda \leq \varepsilon_m$, 2) $\Lambda \leq y_j^k$, and 3) $\Lambda \geq \varepsilon_m + y_j^k - 1$. However, this introduces a computational burden involving $2\times U\times C$ decision variables and additional constraints, rendering it computationally intensive.

Another way is inspired by the branch and bound concept, and we can freeze one variable (constraint linearization). If $y_j^k$ is branched, $2^{C^s}$ sub-problems will be created, the number of which grows exponentially with the number of shared customers. If $\varepsilon_m$ is branched, only $U$ (number of potential meet points) sub-problems are created, which scales obviously much better. Therefore, we choose to branch on the finite set of meet points. More intuitively, by fixing the meet point $m$, we convert the original MINLP to a parametrized ($m$) MILP. The meet point locations hence create $m$ independent branches. By sequentially iterating over them, the global optimum of the original MINLP problem can be found, provided a MILP exact solver is being used. Note that each and every subproblem is still NP-hard to solve. With the adopted branch approach, each sub-problem is a collaborative electric vehicle routing problem with a fixed meet-point (m-CoEVRP). Fig. \ref{FIG:2} demonstrates the parallel computing process of the solution approach. In each subproblem, $\varepsilon_m$ is no longer a decision variable, so that can be removed, and the profit-sharing ratio will be $\alpha_{j}^{m_0}$, and the subproblem can then be reformulated as follows:  

\begin{figure*} [htbp]
	\centering
		\includegraphics[scale=.65]{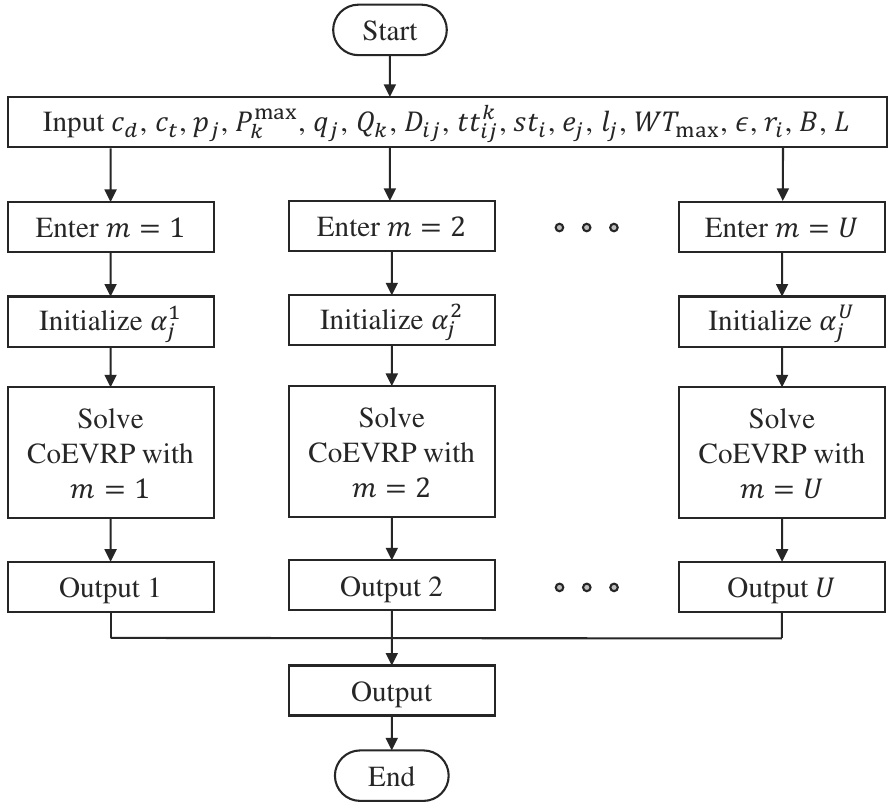}
    \caption{Flowchart of parallel computing of the solution approach}
	\label{FIG:2}
\end{figure*}

\begin{equation}
  \mathrm{min} \;\;\sum_{k\in K}\sum_{i\in N}\sum_{j\in N}c_d{D_{ij}}{x_{ij}^k} + \sum_{k\in K} c_t T^k,
\end{equation}

\noindent\textit{subject to}

\vspace{-\topsep}
\begin{equation} \label{Pknew}
	{\Phi}_k = \sum_{j\in R_k}p_j{y_j^k} + \sum_{j\in R_k}p_j\alpha_j^{m_0}\left(1-y_j^k\right) + \sum_{j\in R \setminus R_k}p_j\left(1-\alpha_j^{m_0}\right)y_j^k - \sum_{i\in N}\sum_{j\in N}c_d{D_{ij}}{x_{ij}^k} - c_t T^k,
\end{equation}

\vspace{-\topsep}
\begin{equation} \label{chargeb1}
	{\delta}_i^k \leq B - b_i^k, \forall{i\in N},{k\in K},
\end{equation}

\vspace{-\topsep}
\begin{equation} \label{chargeb2}
	{\delta}_i^k \leq Bz_i^k, \forall{i\in N},{k\in K},
\end{equation}

\vspace{-\topsep}
\begin{equation} \label{MaxWTn}
	 -{WT_{max}}\leq s_{m_0}^1 - s_{m_0}^2 \leq {WT_{max}},
\end{equation}

\vspace{-\topsep}
\begin{equation} \label{ServiceSeqN}
	s_{m_0}^k - \Gamma\left(1-y_i^k\right) \leq s_i^k, \forall{i\in R_k},{k\in K},
\end{equation}

\vspace{-\topsep}
\begin{equation} \label{OneMPn}
	\sum_{i\in N}{x_{im_0}^k} = 1, \forall{k\in K},
\end{equation}

\vspace{-\topsep}
\begin{equation} \label{SameMPn}
	\sum_{i\in N}{x_{im_0}^1} - \sum_{i\in N}{x_{im_0}^2} = 0,
\end{equation}

and Eq. \eqref{NonPk}, \eqref{BatteryLevel}, \eqref{BatteryCons}, \eqref{capacityV}, \eqref{StartServiceTimeN}-\eqref{CusVisitOnce}, \eqref{yandx}-\eqref{decV2}, \eqref{Profit} remain the same.

Specifically, constraints in Eq. \eqref{Profitk} are linearized to Eq. \eqref{Pknew}; the nonlinear constraints in Eq. \eqref{chargeb} can be easily dealt with by dividing them into two linear constraints in Eq. \eqref{chargeb1} and Eq. \eqref{chargeb2}. The absolute term in Eq. \eqref{MaxWT} is linearized to Eq.\eqref{MaxWTn}. With a fixed meet point, Eq. \eqref{ServiceSeq}, Eq. \eqref{OneMP} and Eq. \eqref{SameMP} can also be simplified to constraints in Eq. \eqref{ServiceSeqN}, Eq. \eqref{OneMPn} and Eq. \eqref{SameMPn}.

As shown in Fig. \ref{FIG:2}, the sub-problems are independent so that they can be solved in parallel. By leveraging parallel computing techniques, the computation time can be significantly reduced, even if the number of potential meet points is considerable. This facilitates real-world implementation, given the fact that the number of meet points is usually rather limited due to requirements such as parking space and regulatory permissions in reality. 

\subsubsection{Metaheuristics with linear programming}

Despite the linearization method introduced above, the MILP subproblems (NP-hard) may necessitate heuristics to address scalability issues. This section presents an approximate approach that integrates heuristics with linear programming to solve the subproblems in Fig. \ref{FIG:2}. We integrate mathematical programming with a heuristic framework, which has been successfully applied for solving various VRP variants in the literature \cite{archetti2014,seyfi2022}. In this paper, we design search-based metaheuristics for route optimization and a linear programming (LP) model for charging schedules. To be specific, given an initial solution, the algorithm includes three interactive modules: 1) an Adaptive Large Neighborhood Search (ALNS) module for route planning, 2) an LP module for charging schedule optimization, and 3) a local search module for further route improvements. 

\begin{figure*} [htbp]
	\centering
		\includegraphics[scale=0.7]{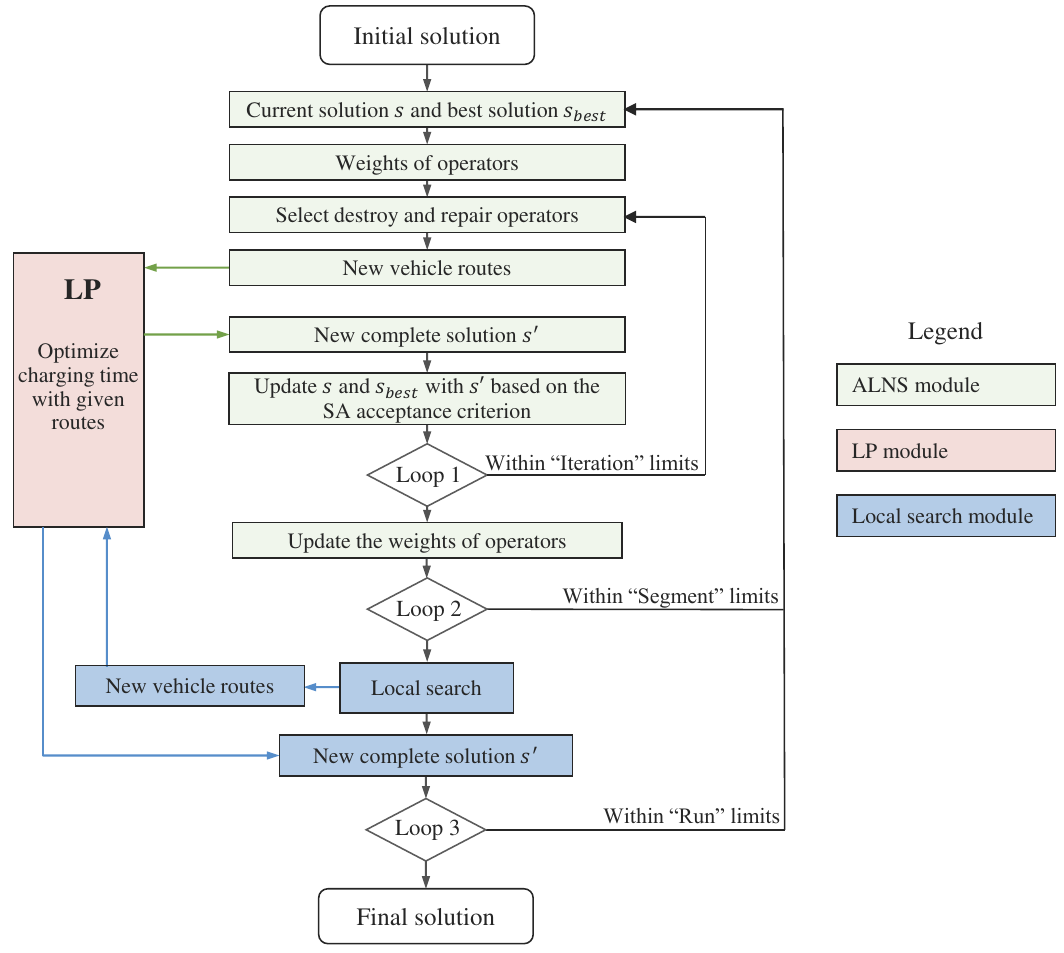}
    \caption{Flowchart of the proposed search-based metaheuristics}
    \label{FIG:4}
\end{figure*}

Before elaborating on the details of the three modules, we first explain the logic of the algorithm framework as illustrated in Fig. \ref{FIG:4}. The modules are utilized in three different layers of loops, namely the "iteration" layer, the "segment" layer, and the "run" layer. Each subsequent layer embeds the previous one, establishing a hierarchical relationship among them. More specifically, each run has a number of segments, and each segment contains a number of iterations.

With an initial solution, the algorithm starts from the bottom layer (iterations) by applying the first part of the ALNS module, i.e., route mutation, which seeks to merely improve EV routes and temporally disregard charging (thus an incomplete solution). Thereafter, within the same loop layer, the LP module is applied to plan the charging schedule, thereby completing the solution. The second part of the ALNS module then evaluates the complete solution and updates operator weights in each segment loop (second layer), where a rewarding mechanism is designed to incentivize better operators. In the outermost layer, the local search module is applied in the run loop to further improve the EV routes, which again needs the LP module to complete the solution. In a nutshell, the CoEVRPMP involves intertwined decision-making in both spatial and temporal domains. The spatial decisions are the sequences of visiting customers and meet points, while the temporal decisions concern the timing and duration of visits to these locations (e.g., charging time). Following this logic, our algorithm divides the solution-finding procedure into the same two domains, with ALNS and local search focusing only on vehicle routes enhancement and the LP module addressing the charging time optimization. Therefore, as shown in Fig. \ref{FIG:4}, whenever a new route is found (either by ALNS or local search), the LP is implemented to complete the solution. The three modules, called in different layers, are designed to iteratively improve the solution in a harmonized and feasibility-guaranteed fashion, which are described as follows.

\textbf{\textit{(i) ALNS module}}

The ALNS module serves as the core of the solving algorithm. ALNS has been widely used and has shown high performance in various VRP variants. We select the ALNS algorithm for its competing performance and flexibility. As demonstrated in recent studies, ALNS could often result in high-quality solutions with acceptable computational run times \cite{keskin2016,hiermann2016,sacramento2019,pelletier2019,chen2021,cheng2023}, which is also the case in our problem (see section \ref{computation} for more details). The flexibility enables us to tailor it to the CoEVRPMP. The original ALNS algorithm was proposed by \citet{ropke2006}, which adopts the principle of removal first and then insertion to find new routes. A set of different removal and insertion operators (destroy and repair operators) are used and assigned with performance-based weights to adaptively improve the solution. More details regarding the standard ALNS algorithm can be found in \citet{ropke2006}.   

We improve the original ALNS algorithm to guarantee solution feasibility, which can be intractable in our problem due to two sets of constraints. First of all, the charging constraints (Eq. \eqref{BatteryCons} to Eq. \eqref{chargeb}) significantly reduce feasible solution space. Secondly, exchanging goods at the meet point puts extra constraints on serving shared and reserved customers, making it even harder to find feasible solutions. To cope with those difficulties, we change the original ALNS algorithm of \citet{ropke2006} in two aspects accordingly. We embed the LP module into the ALNS loops, as shown in Fig. \ref{FIG:4}, so that route finding and charging are handled sequentially, making it easier to find feasible solutions. To resolve the meet-point synchronization, new rules are applied: i) exchanged goods must be delivered after the meet-point, corresponding to Eq. \eqref{ServiceSeq}; ii) to ensure synchronization at the meet-point, we converted the Eq. \eqref{MaxWT} into a penalty function and added it to the objective function. 

\textbf{\textit{(ii) LP module}}
 
As shown in Fig. \ref{FIG:4}, a linear programming (LP) module is applied to optimize charging time whenever a new route solution is obtained. It is worth noting that once the service sequence is determined, for each route, we will be able to calculate the remaining energy at each node ($\widehat{b}_i$) prior to any charging being performed. Therefore, with the service route ($X$) and battery level ($\widehat{b}_i$) information available, the LP module seeks to find optimal charging strategies ($\delta_\ell$) that would minimize the total task time ($s_i$). Note that the $i$ here represents the service sequence instead of the node index, $i\in[1,\zeta]$, where $\zeta$ represents the number of nodes in the route. 
 
After the routing of each vehicle is obtained, the first term of the objective function \eqref{Obj}, the energy consumption cost, is determined. Therefore, the charging battery at each node needs to be optimized to minimize the labor cost, considering the battery and time constraints. The LP model for a single route can thus be formulated as follows:

\vspace{-\topsep}
\begin{equation} \label{LPObj}
	\mathrm{min} \;\;c_t T,
\end{equation}

\textit{subject to}

\vspace{-\topsep}
\begin{equation} \label{LPchargeb1}
	L-\widehat{b}_i \leq \sum_{\ell=1}^i\delta_{\ell} \leq B-\widehat{b}_i, \forall i\in \left[1,\zeta \right],
\end{equation}

\vspace{-\topsep}
\begin{equation} \label{LPservicetime}
	\delta_i + s_i-s_{i+1} \leq -{tt}_{X_iX_{i+1}}-{st}_{X_i}, \forall i\in \left[1,\zeta-1 \right],
\end{equation}

\vspace{-\topsep}
\begin{equation} \label{LPArrivaltime}
	T = s_i, i=\zeta,
\end{equation}

\vspace{-\topsep}
\begin{equation} \label{LPTW}
	e_{X_i} \leq s_i \leq l_{X_i}, \forall i\in \left[1,\zeta \right],
\end{equation}

\vspace{-\topsep}
\begin{equation} \label{LPdelta}
	\delta_i \geq 0, \forall i\in \left[1,\zeta \right].
\end{equation}

In the above model, the objective function Eq. \eqref{LPObj} is a concise version of the original objective function Eq. \eqref{Obj} since the first term is deterministic given a route sequence. The original battery constraints (Eq.\eqref{BatteryCons} and Eq. \eqref{chargeb}) are simplified as Eq. \eqref{LPchargeb1}. The original service time constraints Eq. \eqref{StartServiceTimeN} can be converted to Eq. \eqref{LPservicetime}. The developed LP model has a much lower computational complexity compared to the original MINLP and can be solved by any commercial solver in polynomial time. Without loss of generality, we assume that there are charging facilities at all customer and meet point locations, which also represents the most complex scenario of the studied problem. If charging facilities are not available at some stops, we can easily adapt the model to such simpler cases by restricting charging opportunities defined in Eq. \eqref{LPchargeb1}.

\textbf{\textit{(iii) Local search module}}

The local search module aims to further optimize vehicle routing (and only routes), with the understanding that metaheuristic algorithms can often benefit from extra randomness and disturbances. Through extensive experiments, we discover that the following three operators exhibit the best results: 2-opt \citep{croes1958}, relocate \citep{savelsbergh1992}, and neighbor move. 2-opt and relocate are standard operators, and the neighbor move is specially designed for the studied problem, inspired by the recursive granular algorithm in \citet{moshref2016}. Specifically, the neighbor move is applied to each customer and its pre-determined neighbor customers. Among those customers, one will be selected and relocated as the immediate successor of the ego customer. Note that if an insertion operation is used, it should be guaranteed shared customers should be inserted in their own vehicle's routing or others' routing after the meet point. In contrast, reserved customers can only be inserted in their own vehicle's routing. The local search will accept better new solutions and discard worse solutions.

\section{Numerical experiments} \label{Section4}

In this section, both exact and heuristic-based methods are tested to investigate their viability in various scenarios. To demonstrate the benefits of collaboration, we showcase both small-medium-sized real-world examples and also large-scale problems, using non-collaborative results as benchmarks. We also examine the computational performance of the proposed solution algorithms through a series of numerical experiments. The problem size ranges from 9 customers to 500 customers, representing different use cases.

The experiments are conducted on a standard PC with a six-core Inter(R) Core(TM) i7-8750H CPU at 2.2GHz and 16GB of RAM. The exact method is coded in MATLAB R2021b by using Gurobi 9.5.2 for solving the subproblems. For practical concerns, we impose a limit on the algorithm's runtime, which varies from 0.5h to 100h, depending on the problem size. The heuristic algorithm is also coded in MATLAB R2021b.

\subsection{Real-world case}

In this section, we use a real-world case to demonstrate the merits of collaboration. To cover the vast spectrum of real-world situations, we present comprehensive results with varying vehicle types (EVs or conventional vehicles), time windows, profit thresholds, and numbers of shared customers.

\subsubsection{Case description}

\begin{figure*} [htbp]
	\centering
		\includegraphics[scale=.55]{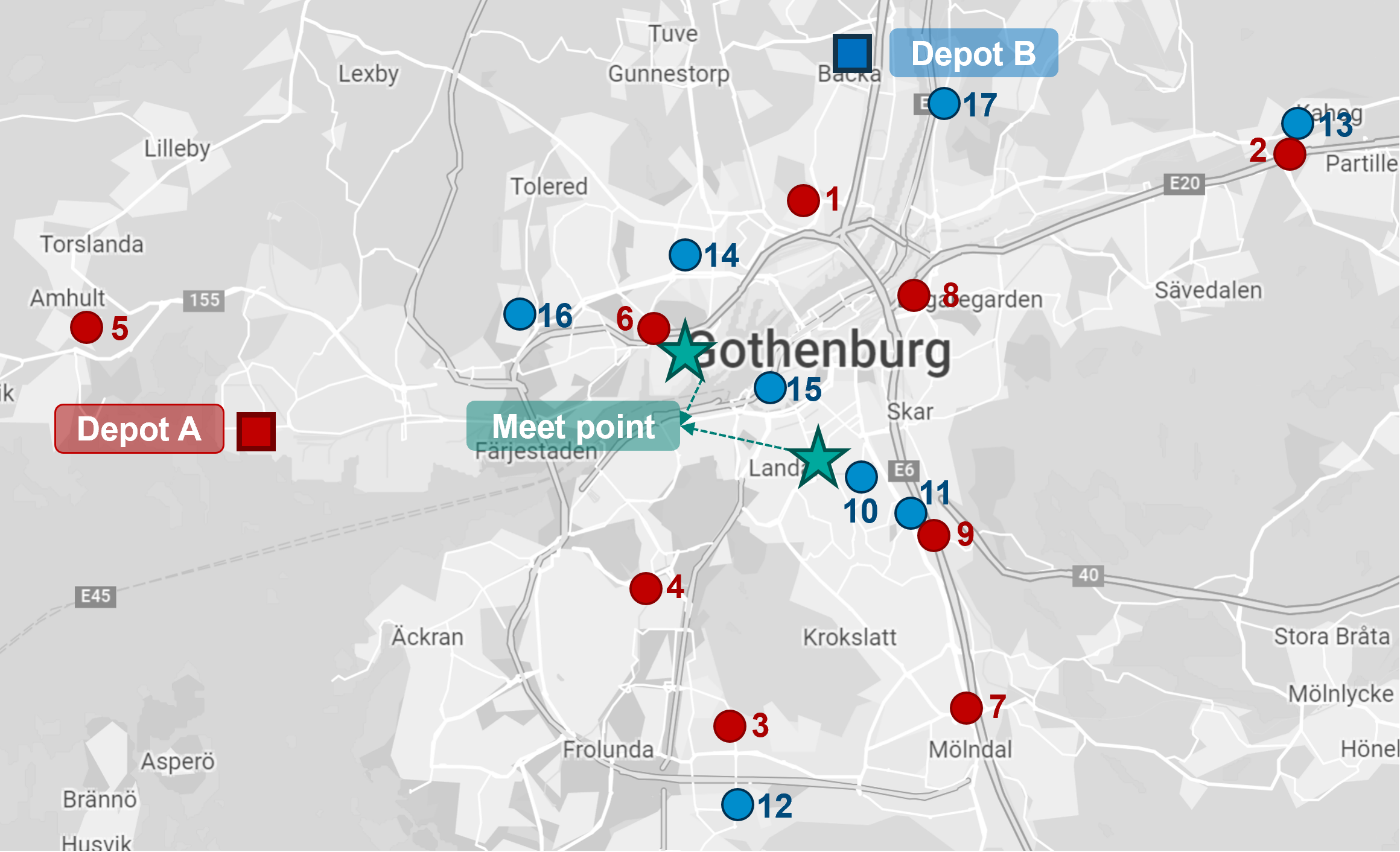}
    \caption{The locations of all nodes}
	\label{FIG:3}
\end{figure*}

The case studies are created based on the real locations of large grocery stores from two companies (namely, ICA and Willys) operating in the city of Gothenburg, Sweden. Both companies routinely deliver goods from depots to their local stores scattered in the city. Fig. \ref{FIG:3} shows the map of interest. Each company has one depot, as marked by the squares. The circles represent local store locations: 9 stores of ICA (red, marked by R) and 8 stores of Willys (blue, marked by B). In our problem, those local stores are the "customers" for the company vehicles to visit. 

For the meet points, we consider places with vacant spaces that can be used for vehicles to meet each other. We assume that the two campuses of the Chalmers University of Technology (Johanneberg and Lindholmen) are optional meet points as marked by the stars in Fig. \ref{FIG:3}. The asymmetric origin-destination distance matrix was obtained from Google Maps API (accessed: 20 April 2022). In addition, the shortest path between any two interest spots is used for distance/time estimation in this study. The distance between nodes is shown in Appendix C (Table \ref{DistanceBetweenNodes}). The time windows related to each store are assumed as those in Table \ref{tbl1}.

\begin{table}[width=.9\linewidth,cols=10,pos=h]
\caption{Time windows of customers (minutes)}\label{tbl1}
\begin{tabular*}{\tblwidth}{@{} LLLLLLLLLL@{} }
\toprule
$R_{\rm{R}}$ & 1 & 2   & 3   & 4   & 5   & 6  & 7 & 8 & 9  \\
$\left[e_j,l_j\right]$ & {[}0,90{]} & {[}30,60{]} & {[}0,90{]}   & {[}30,120{]} & {[}30,120{]} & {[}60,150{]} & {[}60,150{]} & {[}90,180{]} & {[}90,180{]} \\
\midrule
$R_{\rm{B}}$  & 10  & 11 & 12  & 13  & 14  & 15 & 16  & 17  &  \\
$\left[e_j,l_j\right]$ & {[}0,90{]} & {[}0,90{]}  & {[}30,120{]} & {[}60,90{]}  & {[}30,120{]} & {[}60,150{]} & {[}60,150{]} &{[}90,180{]} \\
\bottomrule
\end{tabular*}
\end{table}

Values of other context parameters are determined based on our local survey. Specifically, the amount of service fee ($p_j$) paid by each customer $j$ is 150 in the currency of Swedish Kronor (SEK for short). The unit energy consumption costs ($c_d$) for conventional and electric vehicles are 3SEK/km and 6SEK/km, respectively. And the unit driver salary ($c_t$) is 2.05SEK/min. The average speed ($v$) of both companies' vehicles is assumed to be 40km/h. The travel time (${tt}_{ij}^k$) from node $i$ to node $j$ of vehicle $k$ can is calculated as ${tt}_{ij}^k={D}_{ij}/v$. The service time (${st}_i$) at customer $i$ is 2 minutes; the total unloading and loading time (service time) at meet points is 10 minutes; the maximum waiting time (${WT}_{max}$) for the other vehicle at meet points is 5 minutes. The large positive number ($\Gamma$) is set as 100. For the electric vehicles, we assume total battery capacity $B=60$kWh, the minimum battery is $L=12$kWh (the 20\% of the full battery), unit energy consumption $\epsilon=1$Wh/m, and charging rate $r_i=60$kW.

\subsubsection{Collaboration vs non-collaboration}

We note that the proposed methods can easily adapt to conventional internal combustion engine vehicles, which are still prevailing in the market. In scenarios where conventional vehicles are used, Eq. \eqref{BatteryLevel}, \eqref{chargeb}, \eqref{BatteryCons}, and \eqref{calculateST} can be removed. Moreover, the term ${ST}_i^k$ in constraint \eqref{StartServiceTimeN} will be changed to term ${st}_i$ since there is no charging time. In practice, time windows are sometimes not enforced, in which case we can further remove Eq. \eqref{TW}. 

Due to different problem setups, we solve a few variants of the non-collaborative routing problem: the basic VRP, VRPTW, EVRP, and EVRPTW models. Accordingly, we solve their collaborative counterparts: CoVRPMP, CoVRPMP-TW, CoEVRPMP, and CoEVRPMP-TW. In collaborative cases, it is assumed that all customers can be shared. We use the non-collaboration scenarios as the baselines, where each company's costs and profits are calculated separately. The results are summarized in Table \ref{tbl2}.  

\begin{table}[width=.9\linewidth,cols=13,pos=h]
\caption{Results of collaboration and non-collaboration scenario (SEK)}\label{tbl2}
\begin{tabular*}{\tblwidth}{@{}LLLLLLLLLLLLL@{} }
\toprule
    & \multicolumn{3}{l}{Non-collaboration}  & \multicolumn{9}{c}{collaboration}  \\
    & & & &   & \multicolumn{4}{l}{without profit thresholds} & \multicolumn{4}{l}{with profit thresholds}\\
$k$ & Model & TC  & $\Phi$ & Model   & TC  & ↓(\%)  & $\Phi$ & ↑(\%) & TC  & ↓(\%)   & $\Phi$ & ↑(\%) \\
\midrule
R   & \multirow{2}{*}{VRP}    & \multirow{2}{*}{1223.1} & 720.5  & \multirow{2}{*}{CoVRPMP}    & \multirow{2}{*}{1095.5} & \multirow{2}{*}{10.4}   & 903.3  & 25.4   & \multirow{2}{*}{1124.7} & \multirow{2}{*}{8.0}  & 769.9  & 6.9   \\
B   &   &   & 606.4  &  &   &  & 551.3  & \textcolor{red}{-9.1} &  &  & 655.4  & 8.1   \\
R   & \multirow{2}{*}{VRPTW}  & \multirow{2}{*}{1663.9} & 381.0  & \multirow{2}{*}{CoVRPMP-TW}  & \multirow{2}{*}{1267.5} & \multirow{2}{*}{23.8} & 484.1  & 27.1  & \multirow{2}{*}{1267.5} & \multirow{2}{*}{23.8}  & 484.1  & 27.1  \\
B   &   &  & 505.2  &   &  &   & 798.5  & 58.1   &   &   & 798.5  & 58.1  \\
R   & \multirow{2}{*}{EVRP}   & \multirow{2}{*}{905.6}  & 880.9  & \multirow{2}{*}{CoEVRPMP}   & \multirow{2}{*}{719.2}  & \multirow{2}{*}{20.6}  & 1112.1 & 26.2  & \multirow{2}{*}{736.8}  & \multirow{2}{*}{18.6} & 978.7  & 11.1  \\
B   & &   & 763.5  &   &  &  & 718.7  & \textcolor{red}{-5.9} &   &   & 834.5  & 9.3   \\
R   & \multirow{2}{*}{EVRPTW} & \multirow{2}{*}{1277.7} & 577.9  & \multirow{2}{*}{CoEVRPMP-TW} & \multirow{2}{*}{818.8}  & \multirow{2}{*}{35.9} & 716.6  & 24.0       & \multirow{2}{*}{818.8}  & \multirow{2}{*}{35.9}          & 716.6  & 24.0  \\
B   &   &  & 694.4  &  &   &  & 1014.6 & 46.1  &  &  & 1014.6 & 46.1  \\
\bottomrule
\end{tabular*}
\end{table}

As shown in the table \ref{tbl2}, collaboration reduced the total cost by 8\%-36\% compared to non-collaboration scenarios. These benefits were more pronounced when time windows (24\%-36\%) and electric vehicles (19-36\%) were taken into account. The difference is the largest when electric vehicles are bound by time windows, resulting in a cost reduction of 36\%. We now focus on the impacts of profit thresholds, the number of shared customers, and the length of time windows on the results.

\textbf{\textit{Profit threshold}}

In Table \ref{tbl2}, it could also be found that, without restricting the profit threshold, one company may lose profit as a sacrifice for lower total costs for the two companies. This will, of course, compromise collaboration in real life. For example, in the comparison between VRP and CoVRPMP, the usage of profit thresholds increases total cost but ensures a win-win situation. As shown in the last column of Table \ref{tbl2}, the profits of companies increase by 7\%-58\%, which can lead to a higher willingness to collaborate.

In Assumption \textbf{A3} (refer to Section \ref{Section3}), we highlighted that each company is responsible for setting its own profit thresholds. For the purposes of this paper, we've chosen the non-collaborative profit as our threshold. If a company raises this threshold, it will likely reduce the space for collaboration. On the other hand, by lowering the threshold, there might be more opportunities for collaboration. However, this could come at the expense of individual company profits. At its core, it is about finding a delicate balance between individual gains and collective collaboration for mutual benefit.

\textbf{\textit{Shared customers}}

For the CoEVRPMP problem, different numbers of shared customers are considered. Here, three scenarios are studied: i) 2 shared customers, where only customers 2 and 13 (3 and 12) could be shared, $R^s=\left\{2,13\right\}$ ($R^s=\left\{3,12\right\}$); ii) 4 shared customers, where customers 2, 3, 12, and 13 could be shared, $R^s=\left\{2,3,12,13\right\}$; iii) all customers shared, $R^s=R$, and there are no reserved customers, $R^r=\varnothing$. Taking the result of EVRPTW of non-collaboration as the baseline, the results of different numbers of shared customers are shown in Table \ref{tbl3}.

\begin{table}[width=.9\linewidth,cols=12,pos=h]
\caption{Results with different numbers of shared customers (SEK)}\label{tbl3}
\begin{tabular*}{\tblwidth}{@{} LLLLLLLLLLLL@{} }
\toprule
\multicolumn{2}{l}{\multirow{2}{*}{}}  & \multicolumn{2}{l}{Non-coooperation} & \multicolumn{8}{l}{collaboration}  \\
\multicolumn{2}{l}{}  &  &  & \multicolumn{4}{l}{without profit thresholds}   & \multicolumn{4}{l}{with profit thresholds}\\
$k$ & shared customers $R_s$ & TC & $\Phi$   & TC  & ↓(\%)   & $\Phi$ & ↑(\%) & TC  & ↓(\%)  & $\Phi$ & ↑(\%) \\ \hline
R   & \multirow{2}{*}{$R^s=\left\{2,13\right\}$}  & \multirow{8}{*}{1277.7}   & 577.9 & \multirow{2}{*}{1158.65} & \multirow{2}{*}{9.3}  & 700.0  & 21.1  & \multirow{2}{*}{1245.6} & \multirow{2}{*}{2.51} & 600.3  & 3.9   \\
B   &  &  & 694.4    & &  & 691.3  & \textcolor{red}{-0.4}  &    &  & 704.1  & 1.4   \\
R   & \multirow{2}{*}{$R^s=\left\{3,12\right\}$}   &   & 577.9    & \multirow{2}{*}{1330.9}  & \multirow{2}{*}{\textcolor{red}{-4.2}} & 554.9  & \textcolor{red}{-4.0}  & \multirow{2}{*}{-}      & \multirow{2}{*}{-}    & -      & -     \\
B   &  &  & 694.4    & &  & 664.2  & \textcolor{red}{-4.4}  &    &   & -      & -     \\
R   & \multirow{2}{*}{$R^s=\left\{2,3,12,13\right\}$} &  & 577.9    & \multirow{2}{*}{1058.43} & \multirow{2}{*}{17.2} & 779.9  & 35.0  & \multirow{2}{*}{1058.4} & \multirow{2}{*}{17.2} & 779.9  & 35.0  \\
B   &  &  & 694.4    &  &  & 711.7  & 2.5   &   &  & 711.7  & 2.5   \\
R   & \multirow{2}{*}{$R^s=R$} &  & 577.9    & \multirow{2}{*}{818.8}   & \multirow{2}{*}{35.9} & 716.6  & 24.0  & \multirow{2}{*}{818.8}  & \multirow{2}{*}{35.9} & 716.6  & 24.0  \\
B   &  &   & 694.4    &  &   & 1014.6 & 46.1  &   &   & 1014.6 & 46.1  \\
\bottomrule
\end{tabular*}
\end{table}

In Table \ref{tbl3}, it could be found that, generally, the more shared customers, the better the results in terms of the total cost. With all customers shared, $R^s=R$, the best solution is achieved. However, collaboration may be worse than non-collaboration if only a few customers are shared (see the case $R^s=\left\{3,12\right\}$). This is because the meet point introduces additional travel distances that cannot be compensated by the profits of limited shared customers. In this case, we cannot even find a solution when profit thresholds are applied.

\textbf{\textit{Time windows}}

Last but not least, the impact of time windows is studied. Without losing generality, we examine different combinations of earliest service time $e_i$, length of time windows $\tau_i$, and the total range of time windows of all customers ($\tau=\mathop{max}\limits_i l_i-\mathop{min}\limits_i e_i$, where $l_i$ denotes the latest service time and $\tau_i=l_i-e_i$). For each instance, we randomly generate the earliest service time $e_i$ for each customer, and r* describes the set of $e_i$ for all customers, r*$=:\{e_1, e_2, ..., e_{17}\}$. We can thereby characterize instances as r*-$\tau_i$-$\tau$. Taking r1-60-180 as an example, all customers have the same length of time windows (60 minutes) in this instance, but each has a different $e_i$, and the companies must complete their tasks in 180 minutes. For each combination of $\tau_i$ and $\tau$, we reshuffle r* for five times, and hence r*$ \in \{\rm{r1, r2, r3, r4, r5}\}$. This results in 15 different instances, as shown in Table \ref{tbl4}. 

\begin{table}[width=.9\linewidth,cols=10,pos=h]
\caption{Results with different time window lengths (SEK)}\label{tbl4}
\begin{tabular*}{\tblwidth}{@{} LLLLLLLLLL@{} }
\toprule
           & \multicolumn{3}{l}{Non-collaboration} & \multicolumn{6}{l}{collaboration}                             \\
Instances       & TC        & $\Phi_A$   & $\Phi_B$   & TC    & ↓(\%)         & $\Phi_A$ & ↑(\%) & $\Phi_B$ & ↑(\%) \\
\midrule
r1-60-180  & -         & -          & 638.4      & 940.7  & -    & 876.6    & -     & 732.7    & 14.8  \\
r1-90-210  & 1225.8    & 672.9      & 651.3      & 929.6  & 24.2 & 771.0    & 14.6  & 849.4    & 30.4  \\
r1-120-240 & 1136.9    & 759.7      & 653.4      & 929.6  & 18.2 & 771.0    & 1.5   & 894.4    & 30.0  \\
r2-60-180  & 1109.0    & 733.1      & 707.9      & 878.1  & 20.8 & 929.1    & 26.7  & 742.8    & 4.9   \\
r2-90-210  & 1067.7    & 763.0      & 719.3      & 868.7  & 18.6 & 957.5    & 25.5  & 723.8    & 0.6   \\
r2-120-240 & 1067.7    & 763.0      & 719.3      & 868.7  & 18.6 & 957.5    & 25.5  & 723.8    & 0.6   \\
r3-60-180  & -         & -          & 687.9      & 949.6  & -    & 830.6    & -     & 769.8    & 11.9  \\
r3-90-210  & 1145.3    & 716.0      & 688.7      & 925.8  & 19.2 & 830.6    & 16.0  & 793.6    & 15.2  \\
r3-120-240 & 1064.7    & 796.6      & 688.7      & 919.2  & 13.7 & 889.7    & 11.7  & 741.1    & 7.6   \\
r4-60-180  & -         & -          & 670.0      & 924.5  & -    & 916.0    & -     & 709.5    & 5.9   \\
r4-90-210  & 1225.8    & 630.3      & 693.9      & 922.1  & 24.8 & 904.4    & 43.5  & 723.5    & 4.3   \\
r4-120-240 & 1100.3    & 755.9      & 693.9      & 922.1  & 16.2 & 904.4    & 19.7  & 723.5    & 4.3   \\
r5-60-180  & -         & -          & 618.5      & 1092.2 & -    & 712.1    & -     & 745.7    & 20.6  \\
r5-90-210  & 1188.7    & 729.8      & 631.5      & 937.4  & 21.1 & 850.3    & 16.5  & 762.3    & 20.7  \\
r5-120-240 & 1116.8    & 801.7      & 631.5      & 893.3  & 20.0 & 978.6    & 22.1  & 678.0    & 7.4  \\
\bottomrule
\end{tabular*}
\end{table}

In Table \ref{tbl4}, it can be observed that collaboration is particularly beneficial when dealing with tighter customer time windows. Collaboration can provide feasible solutions even when non-collaboration scenarios fail due to short time windows ($\tau_i=60$). Interestingly, the benefits of collaboration, as reflected by total cost, diminish as the length of time windows increases. This suggests that collaboration could enable more precise and efficient delivery schedules, which not only improves service reliability but also enhances customer satisfaction.

\subsection{Large-scale case studies}

\begin{figure*} [htbp]
	\centering
        \includegraphics[scale=.5]{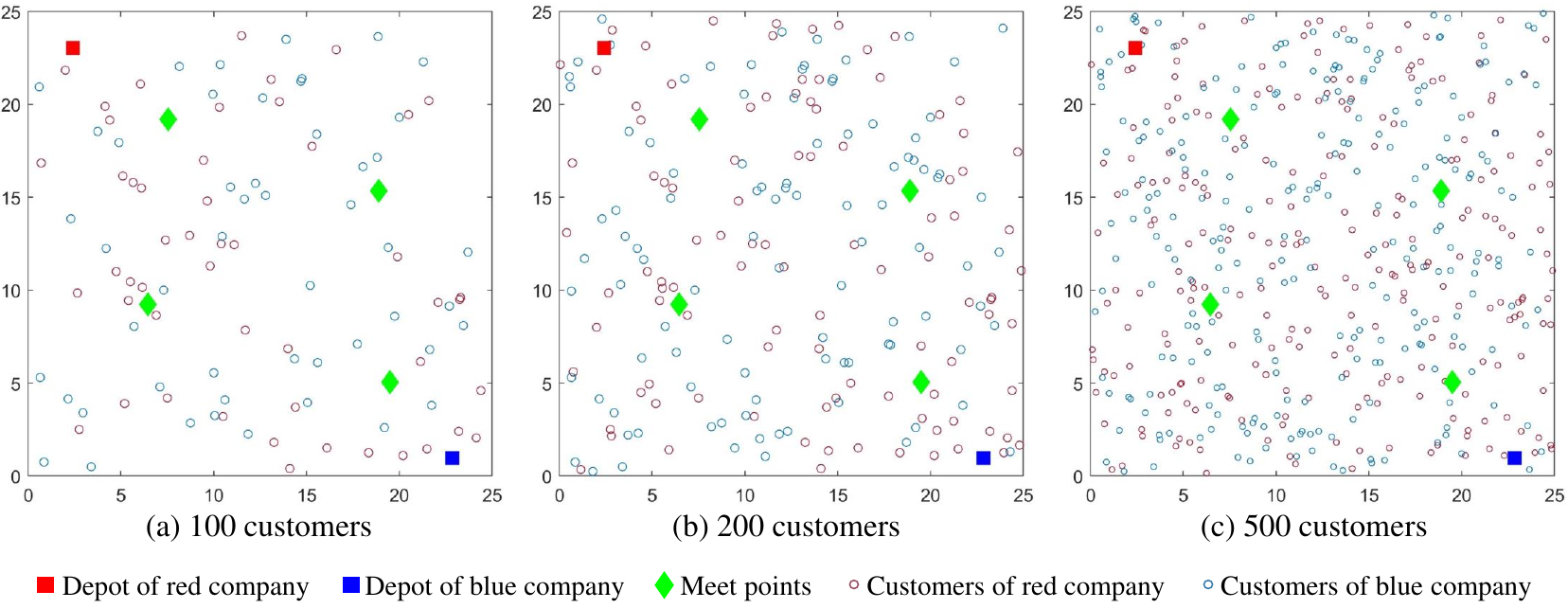}
        \caption{Location map of large-scale cases}
	\label{FIG:6}
\end{figure*}

In this section, we focus on large-scale problems that can only be addressed by the metaheuristics method in practice. We once again consider two companies: red (R) and blue (B). Fig. \ref{FIG:6} illustrates the locations of customers, depots, and meet points with three different problem sizes. Following the same symbolic system as in Fig. \ref{FIG:3}, depots are denoted as squares; diamonds represent meet points, and customers are circles. We locate the depots at polar positions and randomly generate meet points in central zones. Locations of customers are randomly generated, with 100 customers (50 customers for each company), 200 customers (100 customers for each), and 500 customers (250 customers for each) in a 25km $\times$ 25km region. We retain the problem setups as the real-world case, except for 1) the service fee ($p_j$) paid by customer $j$ is 50SEK, 2) the total battery capacity is $B=200$kWh, and 3) the minimum battery is 20\% of the entire battery, $L=12$kWh. To fully examine the impact of collaboration, we solve the problems considering both without TWs and with TWs in those cases. To uphold transparency and facilitate clear communication, we adhere to one vehicle per company. For information on utilizing multiple vehicles, please refer to Appendix A.

\textbf{\textit{Scenarios with time windows}}

In the context of urban logistics, it is often the companies that give optional serving times for customers to choose from (such as DHL and Amazon) instead of the other way around. However, time windows can still vary significantly from case to case, with enormous variants, especially when the number of customers is large. It is thus not practical to examine all possibilities in one study. In this section, we only investigate a non-overlapping two-slot TW\footnotemark[1] set up to showcase the benefits of collaboration in large-scale problems. In the real world, this can represent, for example, a choice of morning vs. afternoon delivery or daytime vs. nighttime delivery. 

\footnotetext[1]{The length of TWs depends on the number of customers, 4, 7, and 14 hours for 100, 200, and 500 customers, respectively.}

\begin{table}[width=.9\linewidth,cols=10,pos=h]
\caption{Results of collaboration and non-collaboration scenario for virtual cases with TWs (SEK)}\label{tbl7}
\begin{tabular*}{\tblwidth}{@{}LLLLLLLLLL@{} }
\toprule
 &      & \multicolumn{3}{l}{Non-collaboration}      & \multicolumn{5}{l}{Collaboration} \\
No. customers   & k    & Model & TC  & $\Phi$ & Model & TC  & ↓(\%)  & $\Phi$ & ↑(\%) \\
\midrule
\multirow{2}{*}{100}  & R  & \multirow{2}{*}{EVRP-TW} & \multirow{2}{*}{3168.5} & 939.8  & \multirow{2}{*}{CoEVRPMP-TW} & \multirow{2}{*}{2440.5} & \multirow{2}{*}{23.0} & 1023.7 & 8.9   \\
& B &  &  & 891.7  &   &    &   & 1535.8 & 72.2  \\
\multirow{2}{*}{200} & R  & \multirow{2}{*}{EVRP-TW} & \multirow{2}{*}{5042.5} & 2448.9 & \multirow{2}{*}{CoEVRPMP-TW} & \multirow{2}{*}{3983.2} & \multirow{2}{*}{21.0} & 2865.0 & 17.0  \\
 & B &    &    & 2508.6 &    &    &   & 3151.9 & 25.6  \\
\multirow{2}{*}{500} & R  & \multirow{2}{*}{EVRP-TW} & \multirow{2}{*}{9158.6} & 7750.0 & \multirow{2}{*}{CoEVRPMP-TW} & \multirow{2}{*}{7691.1} & \multirow{2}{*}{16.0} & 9103.0 & 16.2  \\
& B &  &   & 8012.8 &    &    &    & 8205.9 & 2.4  \\
\bottomrule
\end{tabular*}
\end{table}

\begin{figure*} [htbp]
	\centering
        \includegraphics[scale=.55]{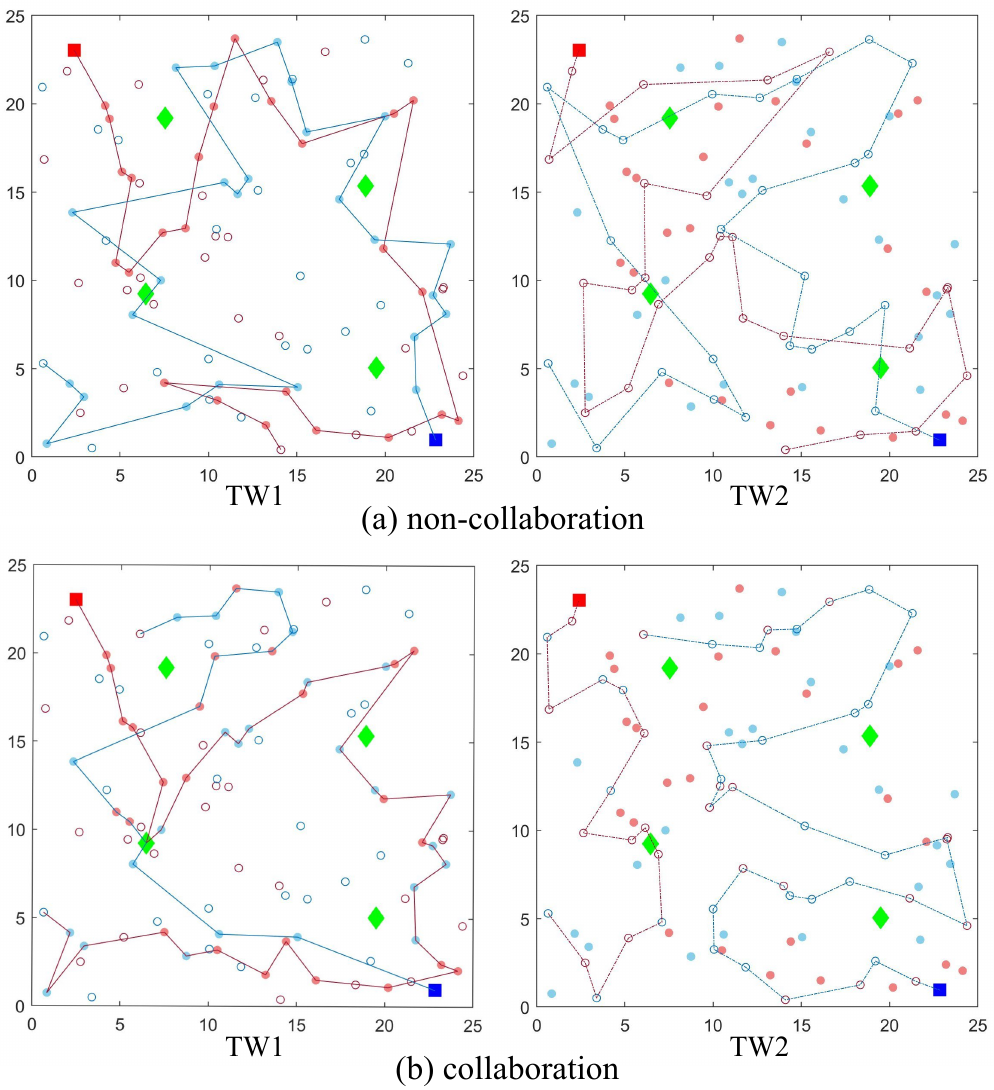}
        \caption{Results of 100 customers with TWs}
	\label{FIG:10}
\end{figure*}

\begin{figure*} [htbp]
	\centering
        \includegraphics[scale=.55]{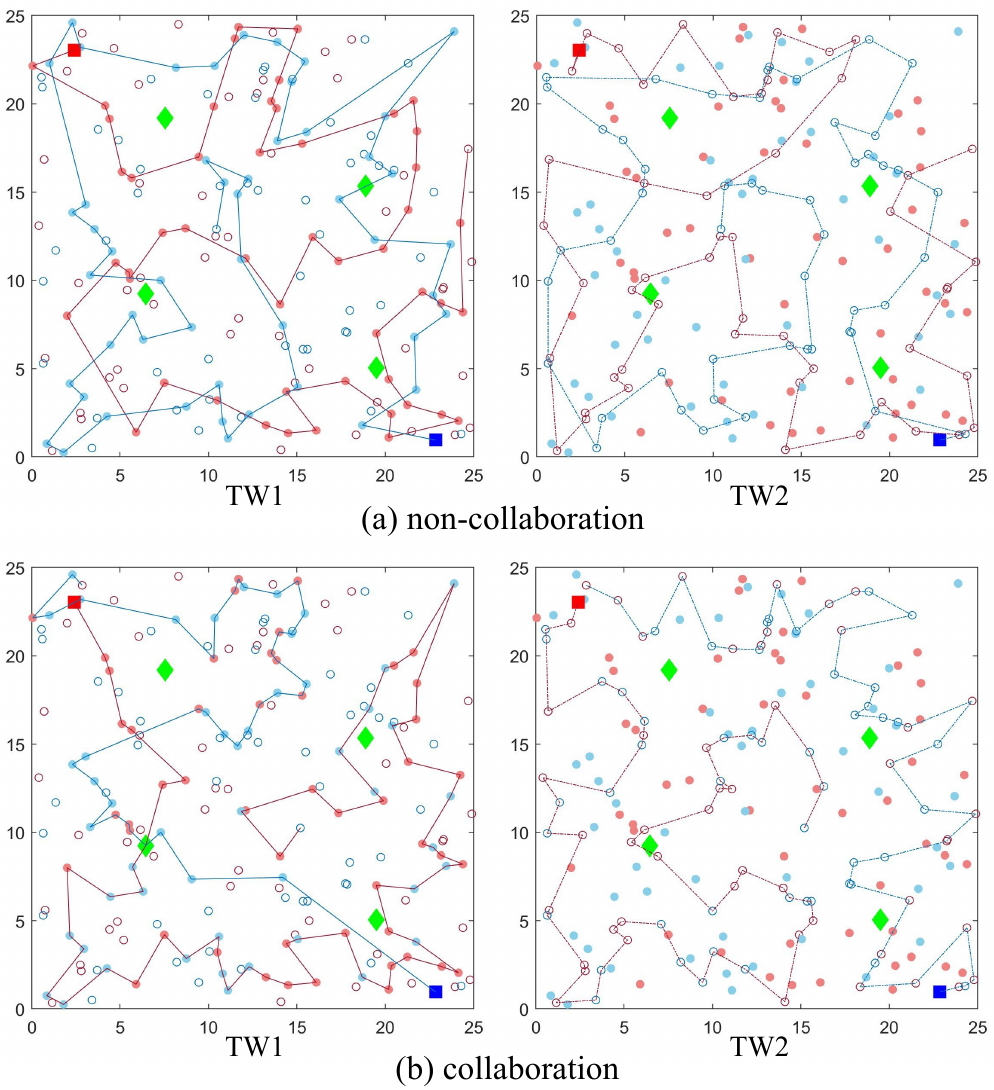}
        \caption{Results of 200 customers with TWs}
	\label{FIG:11}
\end{figure*}

\begin{figure*} [htbp]
	\centering
        \includegraphics[scale=.55]{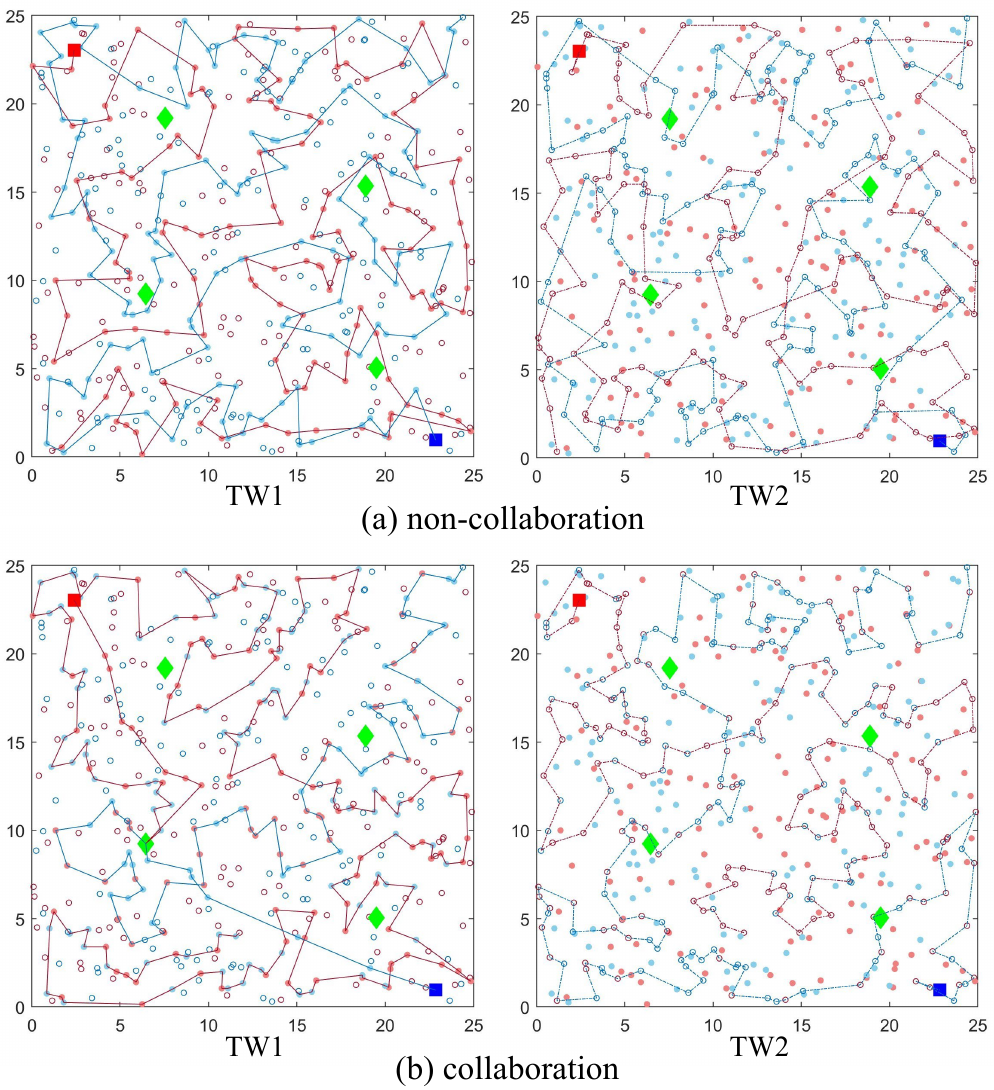}
        \caption{Results of 500 customers with TWs}
	\label{FIG:12}
\end{figure*}

The results are shown in Table \ref{tbl7}. It is evident that collaboration results in significant reductions in the total costs, ranging from 16\% to 23\%. This leads to profit increases for both companies that could incentivize them to collaborate. We could also notice that as the density of customers increases, the benefit of collaboration vanishes. One possible explanation is that it is more rewarding for a company to serve alone if the average energy cost (pure distance-based) between customers is relatively small. In this case, the company can constantly collect profits without traveling too much. On the contrary, if customers are distant from each other, it is advantageous to re-assign the tasks through collaboration so that each company can serve condensed customer areas. This can also be reflected by the routes of vehicles, as shown in Fig. \ref{FIG:10}, Fig. \ref{FIG:11}, and Fig. \ref{FIG:12}. In these figures, the filled dots represent customers who choose the first TW slot, and empty circles denote customers who wish to be served in the later slot. Vehicle routes in these figures indicate that the collaboration tends to separate the pool of customers into two relatively separated clusters for the two companies so that each could focus on a smaller service zone.

Another interesting observation is that, in non-collaboration, vehicle routes often intersect with each other. This means that a vehicle is often bypassing the other company's customers even if they are close and desire the same service time window. In collaborative cases, this issue is resolved by strategically sharing customers based on our models.

\textbf{\textit{Scenarios without time windows}}

\begin{figure*} [htbp]
	\centering
        \includegraphics[scale=.55]{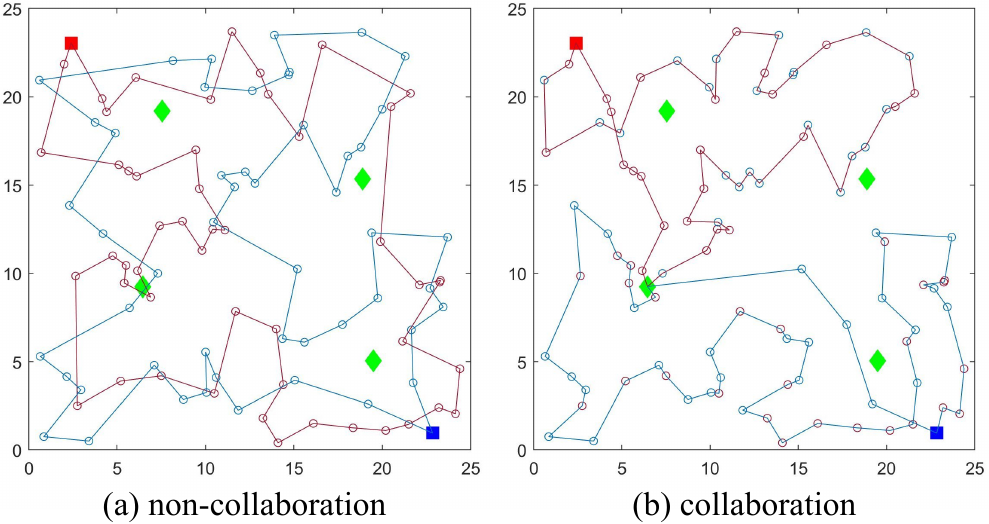}
        \caption{Results of 100 customers without TWs}
	\label{FIG:7}
\end{figure*}

\begin{figure*} [htbp]
	\centering
        \includegraphics[scale=.55]{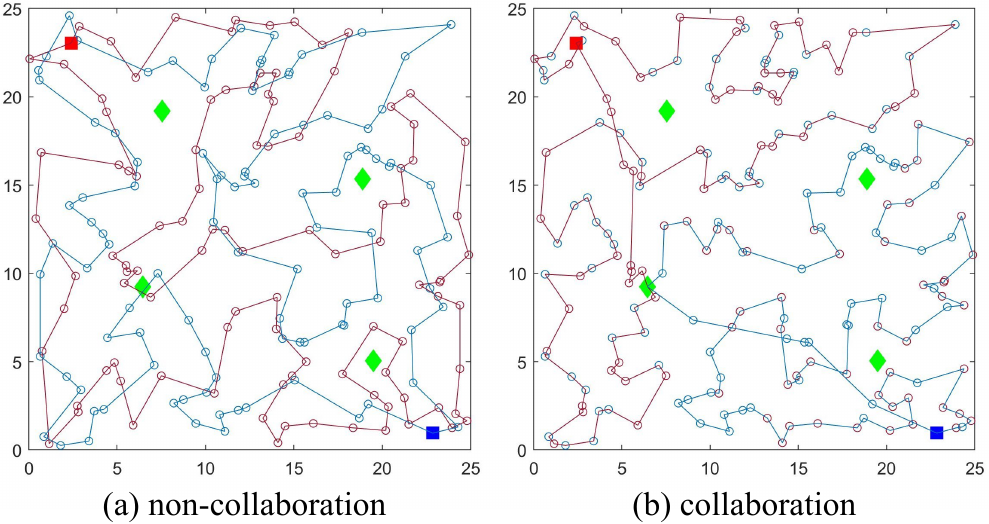}
        \caption{Results of 200 customers without TWs}
	\label{FIG:8}
\end{figure*}

\begin{figure*} [htbp]
	\centering
        \includegraphics[scale=.55]{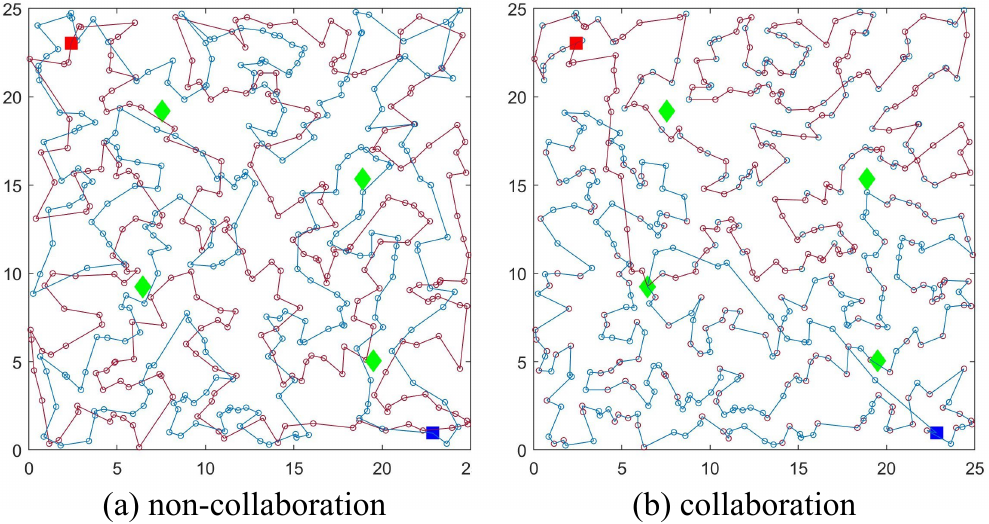}
        \caption{Results of 500 customers without TWs}
	\label{FIG:9}
\end{figure*}

In some scenarios, time windows are not binding or too large to come into effect. We present the results of such cases, as shown in table \ref{tbl6}. The vehicle routes are illustrated in Fig. \ref{FIG:7}, Fig. \ref{FIG:8}, and Fig. \ref{FIG:9}. Compared to cases with TWs, the profits of both companies increase in every scenario since they are freer to plan vehicle routes. In addition, collaboration leads to a clearer separation of service zones without binding time windows.

\begin{table}[width=.9\linewidth,cols=10,pos=h]
\caption{Results of collaboration and non-collaboration scenario for virtual cases without TWs (SEK)}\label{tbl6}
\begin{tabular*}{\tblwidth}{@{}LLLLLLLLLL@{} }
\toprule
 &      & \multicolumn{3}{l}{Non-collaboration}      & \multicolumn{5}{l}{collaboration} \\
No. customers   & k    & Model & TC  & $\Phi$ & Model & TC  & ↓(\%)  & $\Phi$ & ↑(\%) \\
\midrule
\multirow{2}{*}{100} & R  & \multirow{2}{*}{EVRP}    & \multirow{2}{*}{2121.2} & 1477.0 & \multirow{2}{*}{CoEVRPMP}    & \multirow{2}{*}{1719.3} & \multirow{2}{*}{18.9} & 1665.2 & 12.7  \\
   & B &   &   & 1401.8 &   &   &   & 1615.5 & 15.3  \\
\multirow{2}{*}{200} & R  & \multirow{2}{*}{EVRP}    & \multirow{2}{*}{3361.0} & 3287.8 & \multirow{2}{*}{CoEVRPMP}    & \multirow{2}{*}{2674.8} & \multirow{2}{*}{20.4} & 3428.5 & 4.3   \\
& B &   &   & 3351.3 &   &   &   & 3896.6 & 16.3  \\
\multirow{2}{*}{500} & R  & \multirow{2}{*}{EVRP}    & \multirow{2}{*}{6281.1} & 9308.2 & \multirow{2}{*}{CoEVRPMP}    & \multirow{2}{*}{5092.6} & \multirow{2}{*}{18.9} & 9931.6 & 6.7   \\
 & B &   &   & 9410.8 &   &   &   & 9975.8 & 6.0   \\
\bottomrule
\end{tabular*}
\end{table}

In Table \ref{tbl6}, it is obvious that the implementation of collaboration results in approximately 20\% savings in the total costs. It can be clearly seen from Fig. \ref{FIG:7}, \ref{FIG:8} and \ref{FIG:9} that the service areas of each company shrink and the travel distance decreases significantly with collaboration. With collaboration, the company vehicles only need to serve about half of the areas each instead of the entire area. More specifically, when serving 100 customers (see Fig. \ref{FIG:7}), the red vehicle company serves the upper area, while the blue company vehicle serves the lower area. When serving 200 customers (see Fig. \ref{FIG:8}), the red vehicle serves the left top area while the blue vehicle serves the right bottom area in the collaboration scenario. With the increase of the number of customers to 500 (see Fig. \ref{FIG:9}), we receive similar to the 100 customer case pattern.

\subsection{Computational performance} \label{computation}

In this subsection, we examine the computational performance of the two solution approaches. Based on our experiments, the removal fraction $\rho$ of 0.3 renders the best performance, and other parameters are tuned based on the method proposed in \citet{ropke2006}. We use the following naming format $R$-$R_1^r$-$R_1^s$-$R_2^r$-$R_2^k$ to denote different instances. Taking instance 9-0-5-0-4-6 as an example, there are 9 customers in total, including 0 reserved and 5 shared customers of company 1, 0 reserved and 4 shared customers of company 2 as well. The key performance metrics are defined and explained in Table \ref{tbl10}. Table \ref{tbl11} compares the performance of the proposed metaheuristic algorithm with the exact method.

\begin{table}[width=.9\linewidth,cols=2,pos=h]
\caption{Abbreviation of experiment indicators and definition}\label{tbl10}
\begin{tabular*}{\tblwidth}{@{} LL@{} }
\toprule
Abbreviation & Definition\\
\midrule
$S_{Exact}$ & The best feasible objective value found by Gurobi solver in a preset running time\\
$S_{Meta}$ & The best feasible objective value found by the metaheuristic after a preset number of iterations\\
${Imp}_{Exact-Meta}$ & The improvement of ${S}_{Meta}$ compare to ${S}_{Exact}$, which is calculated by: $\frac{S_{Exact}-S_{Meta}}{S_{Exact}}$\\
${CPU}_{Meta}$ & The computation duration of the metaheuristic\\
${CPU}_{Exact}$ & CPU time for solving the MILP model by Gurobi\\
\bottomrule
\end{tabular*}
\end{table}

\begin{table}[width=.9\linewidth,cols=6,pos=h]
\caption{Computational results for instances}\label{tbl11}
\begin{threeparttable}
\begin{tabular*}{\tblwidth}{@{} LLLLLL@{} }
\toprule
Instances & ${S}_{Exact}(\rm{SEK})$ & ${S}_{Meta}(\rm{SEK})$ & ${Imp}_{Meta-Exact}(\%)$ & ${CPU}_{Exact}\rm{(s)}$ & ${CPU}_{Meta}\rm{(s)}$ \\
\midrule
9-0-5-0-4 & \textbf{611.2} & \textbf{611.2} & 0 & 13.4 & 3.7 \\
10-0-5-0-5 & \textbf{512.6} & \textbf{512.6} & 0 & 157.0 & 5.3 \\
10-2-3-2-3 & \textbf{585.7} & \textbf{585.7} & 0 & 38.2 & 7.2 \\
10-4-1-4-1 & \textbf{609.5} & \textbf{609.5} & 0 & 16.5 & 10.9 \\
15-0-7-0-8 & 656.0 & 656.0 & 0 & 5,400 & 16.1 \\
15-2-5-2-6 & 666.7 & 667.9 & \textcolor{red}{-0.2} & 5,400 & 18.9 \\
15-3-4-4-4 & 690.3 & 691.6 & \textcolor{red}{-0.2} & 5,400 & 32.2\\
20-0-7-0-13 & 720.2 & 707.6 & 1.8 & 10,800 & 19.7 \\
20-2-5-5-8 & 781.7 & 742.5 & 5.3 & 10,800 & 46.1\\
20-6-1-11-2 & 835.3 & 778.79 & 7.3 & 10,800 & 84.7\\
30-0-15-0-15 & 1008.0 & 1007.6 & 0.0 & 14,400 & 34.6\\
30-8-7-8-7 & 1333.8 & 1207.9 & 10.4 & 14,400 & 77.9\\
40-0-20-0-20 & 1232.6 & 1228.6 & 0.3 & 36,000 & 108.3\\
40-10-10-10-10 & 1626.7 & 1591.0 & 2.2 & 36,000 & 305.4\\
50-0-25-0-25 & - & 1403.4 & - & 43,200 & 160.2 \\
60-0-30-0-30 & - & 1354.6 & - & 54,000 & 227.8 \\
60-20-10-20-10 & - & 1539.5 & - & 54,000 & 261.5 \\
80-0-40-0-40 & - & 1559.4 & - & 72,000 &  389.1\\
80-20-20-20-20 & - & 1901.0 & - & 72,000 & 456.3 \\
100-0-50-0-50 & - & 1747.3 & - & 108,000 & 529.2 \\
100-25-25-25-25 & - & 1964.8 & - & 108,000 &  731.9 \\
200-0-100-0-100 & - & 2706.2 & - & 216,000 & 1137.8 \\
500-0-250-0-250 & - & 5150.0 & - & 360,000 & 2898.6 \\
\bottomrule
\end{tabular*}
\begin{tablenotes}
\item The proven optima solutions are indicated in boldface. The negative improvement percentages in column 4, which means the solutions of metaheuristics are worse, are marked in red.
\end{tablenotes}
\end{threeparttable}
\end{table}

As shown in Table \ref{tbl11}, the metaheuristic algorithm generally outperforms the exact approach in terms of computational time and solution quality. Both approaches are able to solve the small-scale problem to optimality. As the problem size increases, it becomes increasingly unfeasible to find optimal solutions via the exact method within reasonable time limits. To this end, we apply a relative indicator, ${Imp}_{Meta-Exact}$, to compare the approximate solutions derived from both approaches. It can be found that the near-optimal solutions obtained by the two solution approaches are very close. Specifically, the metaheuristic reveals a modest gap at most, by a marginal 0.2\% from the exact method. Yet, in more cases, it leads to improved solutions, outperforming the exact method by a significant margin of up to 10.4\%. However, when the problem size exceeds 40, the exact method fails to return any feasible results. In such cases, the heuristic-based approach is the only viable option for practical implementation, which aligns with findings from plenty of existing studies \citep{wang2019,ma2023,xia2023}.   

\section{Conclusion and future work} \label{Section5}

This paper introduced and analyzed a CoEVRPMP. We integrated profit sharing into routing planning with explicit considerations of practical constraints such as charging, time windows, vehicle capacity, and meet-point synchronization. Two solving methods are developed for the formulated CoEVRPMP, i.e., an exact method and a metaheuristic algorithm. The computational performance of the proposed solution methods is examined. Numerical experiments based on real-world cases and large-scale cases are conducted to demonstrate the benefits of collaboration and examine the impacts of profit threshold, the proportion of shared customers, and time windows.

Our results indicate that collaboration with meet points can effectively reduce the total cost by 8\%-36\% compared to the non-collaboration scenarios. In practice, this is reflected in fewer overlaps and intersections between the routes of the two companies. Our results also suggest that in cases where only a few customers are shared, collaboration may not necessarily be beneficial; however, the collaborative benefit increases as the companies share more customers. This indicates significant potential for application in real-world scenarios involving larger customer bases. 

The current work also seeks a win-win collaboration by considering profit threshold constraints for participating companies. Removing the threshold could potentially lead to lower total costs but bears the risk of sacrificing the profit of one party and thus compromising collaboration. Moreover, collaborative routing is more beneficial to electric vehicles than conventional vehicles, especially when time windows are considered.

The results of our experiments also indicate that collaboration significantly outperforms non-collaborative solutions in scenarios with tighter time windows. Here, profit savings were amplified when time windows were shorter compared to longer ones. This hints that the duration of the time window may be effectively reduced in collaborative scenarios, which would, in turn, enhance customer satisfaction since customers usually prefer shorter standby times. In this way, collaborating companies may be able to offer more reliable and precise time windows together. 

In addition, our experiments indicate that the proposed method is most effective in service areas characterized by low demand density. In essence, when customers are dispersed, collaboration can yield greater cost reductions. This observation underscores the need for future research to thoroughly examine the relationship between customer network topology and the benefits derived from collaboration.

Our idea of meet points can be used not only in horizontal collaboration among companies but also in vertical collaboration or collaboration scenarios within a single company. It can be the collaboration within a multimodal system, which means the carrier can be of any type, and goods can be transferred from one mode to another at meet points. For example, this could involve collaboration between electric vehicles and cargo bikes or between trucks and drones. However, if the situation considers only the transfer of goods from one mode to another without exchanging goods, the point of encounter becomes a transshipment node.

We only consider two companies in the main study, each with one vehicle to focus on the core problem of collaboration. The model presented in this paper has robust scalability that can be extended in several different ways. First, we present the potential extension to multiple vehicles in the Appendix. Subsequent research can encompass multiple companies with multiple vehicles in collaborative routing problems involving meet points. In such scenarios, goods can be exchanged at customers' premises, which can be further examined and discussed within the routing models. Second, the parallel branching structure we proposed can be used for not only exact algorithms but also metaheuristics so that computational time can be saved significantly by utilizing parallel computing. Besides, more comprehensive and complex energy consumption estimation methods could be considered in EVRPs and CoEVRPs rather than distance-based ones.

Moreover, the profit increases of the collaborating companies are, at times, uneven in our numerical experiments. Although thresholds ensure a win-win situation, one of the companies may still decline collaboration if the profit increase of the other company is significantly larger. To this end, future research can investigate the balance of profits for the companies in order to make them more willing to collaborate. Last but not least, the problem we studied entails a central authority to coordinate the collaboration. How to enable collaboration without such a trusted consolidator needs further exploration.

\section*{Acknowledgments}
This work was in part supported by the Transport Area of Advanced within Chalmers University of Technology through the project COLLECT: Horizontal Cooperation in Urban Distribution Logistics - a Trusted - Cooperative Electric Vehicle Routing Method.

\clearpage
\appendix

\section*{Appendix A} \label{AppendixA} 

Here, we extend the developed model in Section \ref{Section3} to multiple vehicle scenarios. For company $k$, $v_k$ vehicles are present, with $v_k$ belonging to the set $V_k=\left\{V_1, V_2 \right\}$, which is a subset of the overall vehicle set $V$. The number of vehicles can vary among companies. We assume that exchanges occur exclusively between distinct companies, negating the requirement for all vehicles to visit meet points. Compared with the original model, two additional decisions should be added and considered: 1) pairing two vehicles and determining the meet point for their rendezvous, and 2) determining if a customer transfer is necessary at the meet point and identifying the appropriate location for it.

To accommodate multiple delivery vehicles, we slightly modify the model outlined in Section \ref{Section3.1}. Firstly, we introduce two new decision variables. Secondly, we update the existing decision variables associated with vehicle $k$ to refer to vehicle $v_k$. Thirdly, some variables are removed, while others are retained. Specifically, the decision variable $\chi_m^{v_1v_2}$ signifies whether vehicles $v_1$ and $v_2$ meet at meet point $m$, where $v_1$ and $v_2$ are part of the vehicle set $V$, and $m$ belongs to the meet point set $M$. The decision variable $\varepsilon_j^m$ indicates whether customer $j$ is transferred at meet point $m$. The original decision variables $x_{ij}^k$, $z_i^k$, $s_i^k$, $b_i^k$, $\delta_i^k$, and $ST_i^k$ are updated to $x_{ij}^{v_k}$, $z_i^{v_k}$, $s_i^{v_k}$, $b_i^{v_k}$, $\delta_i^{v_k}$, and $ST_i^{v_k}$, respectively. Furthermore, we retain decision variables $\alpha_j^m$ and remove $y_j^k$ and $\varepsilon_m$.

Undoubtedly, some equations need slight adjustments due to replacing $k$ with $v_k$. For equations \eqref{Obj}, \eqref{BatteryCons}, \eqref{BatteryLevel}, \eqref{chargeb}, \eqref{capacityV}, \eqref{StartServiceTimeN}, \eqref{calculateST}, \eqref{Arrivaltime}, \eqref{TW}, \eqref{CusVisitOnce}, \eqref{depotSandE}, \eqref{depotSandE2}, \eqref{reservedCus}, and \eqref{conservation}, apart from changing $k$ to $v_k$, the formulas stay the same. Equations \eqref{NonPk} and \eqref{Calpha} remain unchanged. Additionally, we present equations that underwent significant modifications below.

\vspace{-\topsep}
\begin{equation} \label{Profitm}
\begin{aligned}
    {\Phi}_k  = &\sum_{j\in R_k}p_j\left(1-\sum_{m \in M}{\varepsilon_j^m}\right) + \sum_{m\in M}\sum_{j\in R_k}p_j\alpha_j^{m}\varepsilon_j^m + \sum_{m\in M}\sum_{j\in R\setminus R_k}p_j\left(1-\alpha_j^{m}\right)\varepsilon_j^m\\
    &- \sum_{v_k \in V_k}\sum_{i\in N}\sum_{j\in N}c_d{D_{ij}}{x_{ij}^{v_k}} - \sum_{v_k \in V_k}c_t T^{v_k},
\end{aligned}
\end{equation}

\vspace{-\topsep}
\begin{equation} \label{CusTransOnce}
    \sum_{m\in M} \varepsilon_j^m \leq 1, \forall{j \in R},
\end{equation}

\vspace{-\topsep}
\begin{equation} \label{OneMPm}
	\sum_{i\in N}\sum_{m\in M}{x_{im}^{v_k}} \leq 1, \forall{v_k\in V_k},{k\in K},
\end{equation}

\vspace{-\topsep}
\begin{equation} \label{OneMPm1}
    \sum_{m\in M} \sum_{v_2\in V_2} \chi_m^{v_1v_2} = \sum_{i\in N}\sum_{m\in M}{x_{im}^{v_1}}, \forall{v_1 \in V_1},
\end{equation}

\vspace{-\topsep}
\begin{equation} \label{OneMPm2}
    \sum_{m\in M} \sum_{v_1\in V_1} \chi_m^{v_1v_2} = \sum_{i\in N}\sum_{m\in M}{x_{im}^{v_2}}, \forall{v_2 \in V_2},
\end{equation}

\vspace{-\topsep}
\begin{equation} \label{SameMPm}
	\left(\sum_{i\in N}{x_{im}^{v_1}} - \sum_{i\in N}{x_{im}^{v_2}}\right)\chi_m^{v_1v_2} = 0, \forall{m\in M},{v_1\in V_1},{v_2\in V_2},
\end{equation} 

\vspace{-\topsep}
\begin{equation} \label{varey}
    \sum_{m \in M}\varepsilon_j^m + \sum_{v_k \in V_k}\sum_{i \in N} x_{ij}^{v_k}= 1, \forall{j \in R_k},{k \in K},
\end{equation}

\vspace{-\topsep}
\begin{equation} \label{varechi}
    \sum_{i\in N}\sum_{j \in R_2}x_{ij}^{v_1} + \sum_{i\in N}\sum_{j \in R_1}x_{ij}^{v_2}\geq  \sum_{m \in M}\chi_m^{v_1v_2}, \forall{v_1\in V_1},{v_2\in V_2}
\end{equation}

\vspace{-\topsep}
\begin{equation} \label{xvarechi}
    \sum_{i \in N}x_{ij}^{v_1}\varepsilon_j^m \leq \sum_{v_2 \in V_2} \chi_m^{v_1v_2}, \forall{m\in M},{j\in R},{v_1\in V_1}
\end{equation}

\vspace{-\topsep}
\begin{equation} \label{xvarechi2}
    \sum_{i \in N}x_{ij}^{v_2}\varepsilon_j^m \leq \sum_{v_1 \in V_1} \chi_m^{v_1v_2}, \forall{m\in M},{j\in R},{v_2\in V_2}
\end{equation}

\vspace{-\topsep}
\begin{equation} \label{MaxWTm}
	\left|s_{m}^{v_1} - s_{m}^{v_2}\right|-\Gamma \left(1-\chi_m^{v_1v_2}\right) \leq {WT_{max}}, \forall{m\in M},{v_1 \in V_1},{v_2 \in V_2}
\end{equation}

\vspace{-\topsep}
\begin{equation} \label{ServiceSeqm}
	s_{m}^{v_k} - \Gamma \left(1-\varepsilon_j^m\right) \leq s_j^{v_k}, \forall{j\in R-R_k},{v_k\in V_k},{k\in K},{m\in M}
\end{equation}

For the constraints, we update them as follows. Eq.\eqref{Profitm} replaces Eq. \eqref{Profitk} to represent the profit of company $k$. Eq. \eqref{CusTransOnce} ensures that each customer can be transferred at most once. Eq. \eqref{OneMPm}-\eqref{OneMPm2} ensure that each vehicle visits at most one meet point. Among them, Eq. \eqref{OneMPm} is to instead Eq. \eqref{OneMP}. Eq. \eqref{SameMPm} guarantees that vehicles $v_1$ and $v_2$ will meet at the same meet point $m$ if they are designated to exchange goods, which has a similar meaning to Eq. \eqref{SameMP}. Eq. \eqref{varey} ensures that if customer $j$ requires a transfer at a meet point, the other company should handle the service; if no transfer is needed, the customer should be served by the original company. Eq. \eqref{varechi} ensures that when vehicles $v_1$ and $v_2$ converge at meet points, at least one customer transfer occurs; otherwise, vehicles do not meet there. Eq. \eqref{xvarechi} and \eqref{xvarechi2} ensure that if vehicles meet at a meet point, then the goods in both vehicles can only be transferred at this point if needed. The waiting time at the meet point needs to be guaranteed within $WT_{max}$ by Eq. \eqref{MaxWTm}, which replaces Eq. \eqref{MaxWT}. Eq. \eqref{ServiceSeqm} replaces Eq. \eqref{ServiceSeq}, ensuring the service sequence that the exchanged goods must be served after the meet points.

\vspace{-\topsep}
\begin{equation} \label{decV1m}
	x_{ij}^{v_k},z_{i}^{v_k}\in \left\{ 0,1 \right\}, \forall{i\in N},{j\in N},{v_k\in V_k}, {k\in K},
\end{equation}

\vspace{-\topsep}
\begin{equation} \label{decV2m}
	s_{i}^{v_k},b_{i}^{v_k},{\delta}_{i}^{v_k},{ST}_{i}^{v_k} \geq 0, \forall{i\in N},{v_k\in V_k}, {k\in K},
\end{equation}

\vspace{-\topsep}
\begin{equation} \label{varepchim}
	\varepsilon_j^m, \chi_m^{v_1v_2}\in \left\{ 0,1 \right\}, \forall{j\in R},{v_1 \in V_1},{v_2 \in V_2},{m\in M}.
\end{equation}

Moreover, Eq. \eqref{Profit}, \eqref{decV1m}, \eqref{decV2m}, and \eqref{varepchim} are the decision variable domains.

\section*{Appendix B} \label{AppendixB}

A table with all notations used in this paper is presented below, including abbreviations, sets, parameters, and decision variables.

\begin{longtable}{ll}
\caption{Mathematical notation}\label{tbnotation} \\
\hline
\textbf{Abbreviation} &   \\
VRP & Vehicle routing problem\\
EVRP & Electric vehicle routing problem\\
TW & Time windows\\
EVRPTW & Electric vehicle routing problem with time windows\\
CoVRP & Collaborative vehicle routing problem\\
CoVRPMP & Collaborative vehicle routing problem with meet points\\
CoVRPMP-TW & Collaborative vehicle routing problem with meet points and time windows\\
CoEVRPMP & Collaborative electric vehicle routing problem with meet points\\
CoEVRPMP-TW & Collaborative electric vehicle routing problem with meet points and time windows\\
m-CoEVRP & Collaborative electric vehicle routing problem with a fixed meet point\\
MDVRP & Multi-depot vehicle routing problem\\
PD & Pickup and delivery\\
PDP & Pickup and delivery problem\\
MILP & Mixed-integer linear programming\\
MINLP & Mixed-integer nonlinear programming\\
ALNS & Adaptive Large Neighborhood Search\\
LP & Linear programming\\
TC & Total cost\\
SEK & Swedish Kronor\\
\hline
\textbf{Sets:} &                    \\
$R_k^{r}$  & The reserved customers of company $k$\\
$R_k^{s}$  & The shared customers of company $k$\\
$R_k$  & The customers of company $k$, $R_k^{r} \cup R_k^{s} = R_k$ \\
$R$  & All the customers\\
$K$  & Vehicles and companies, $k \in K$, in which $k$ is the index of companies/vehicles, $K = \left\{1,2\right\}$\\
$M$  & The meet points, $m\in M$, in which $m$ is the index of meet points\\
$O$  & The depots of companies, $o_k\in O$, $O^+$ and $O^-$ are the start and end depots, $O^+\cup O^- = O$\\
$N$  & All nodes, $N = R\cup M\cup O$\\
\hline
\textbf{Parameters:} &      \\
$c_d$ & Unit energy consumption cost (SEK/km)\\
$c_t$ & Unit driver salary (SEK/min)\\
$D_{ij}$ & Distance from node $i$ to node $j$ (km)\\
$p_j$ & The service fee customer $j$ pays for the delivery service (SEK)\\
$P_k^{min}$ & Minimum profit threshold of company $k$ (SEK)\\
$q_j$ & Demand of customer $j$\\
$\alpha_j^m$ & Profit ratio of customer $j$ exchange goods at meet-point $m$\\
$Q_k$ & Capacity of vehicle $k$\\
${tt}_{ij}^k$ & Travel time from node $i$ to node $j$ for vehicle $k$ (minutes)\\
${st}_i$ & Service time of goods at node $i$ (minutes)\\
$\left[e_i,l_i\right]$ & Time window within which the vehicle should begin to serve node $i$ (minutes)\\
${WT}_{max}$ & Maximum time of the first arrival vehicle at the meet-point waiting for another (minutes)\\
$B$ & Total battery capacity (Wh)\\
$L$ & Minimum battery (Wh)\\
$\epsilon$ & Unit energy consumption per distance (W/km)\\
$r_i$ & Charging rate of charging node $i$ (W)\\
$U$ & Number of potential meet points\\
$C$ & Number of customers\\
$C^s$ & The number of shared customers\\
$C^r$ & The number of reserved customers\\
$\Gamma$ & A large positive number\\
$X$ & The service sequence of a single route\\
$\widehat{b}_i$ & Remaining battery at service sequence position $i$ before any charging performed (Wh)\\
$tt_{X_iX_{i+1}}$ & Travel time from node $X_i$ to node $X_{i+1}$ (minutes)\\
$st_{X_i}$ & Service time of goods at node $X_i$ (minutes)\\
$\left[e_{X_i},l_{X_i}\right]$ & Time window within which the vehicle should begin to serve node $X_i$ (minutes)\\
\hline
\textbf{Decision variables:} &   \\
$x_{ij}^k$ & 1 if vehicle $k$ delivers from node $i$ to node $j$; otherwise 0\\
$y_j^k$ & 1 if customer $j$ is served by vehicle $k$; otherwise 0\\
$z_i^k$ & 1 if vehicle $k$ charges at node $i$, otherwise 0\\
$\varepsilon_m$ & 1 if vehicles choose to meet at meet point $m$\\
$b_i^k$ & The remaining energy in the battery of vehicle $k$ when arriving at node $i$ (Wh)\\
${\delta}_i^k$ & The charging battery of vehicle $k$ at node $i$ (Wh)\\
${ST}_i^k$ & Time for serving goods and charging of vehicle $k$ at node $i$ (minutes)\\
$s_i^k$ & Time at which vehicle $k$ begins service at node $i$ (minutes)\\
$T^k$ & Arrival time of vehicle $k$ at the end depot (minutes)\\
${\Phi}_k$ & Total profit of company $k$ (SEK)\\
$\delta_\ell$ & The charging battery at service sequence position $\ell$\\
$s_i$ & Start service time at service sequence position $i$ (minutes) \\
$T$ & Arrival time at the end depot (minutes)\\
\hline
\end{longtable}

\section*{Appendix C} \label{AppendixC}

Table \ref{DistanceBetweenNodes} shows the actual distance between nodes in the case study, including customer nodes, meet points, and depots. Among them, 1-17 are customer points, 1-9 are customers of company R, and the rest are company B's customers; m1 and m2 are meet points; D1 and D2 are the depots of company R and B, respectively.

\begin{sidewaystable}[!htp]
\small
\caption{The actual distance between nodes (km)} \label{DistanceBetweenNodes}
\centering
\begin{tabular*}{\linewidth}{LLLLLLLLLLLLLLLLLLLLLL}
\toprule
\multicolumn{1}{c}{} & \multicolumn{1}{c}{1} & \multicolumn{1}{c}{2} & 3    & 4    & 5    & 6    & 7    & 8    & 9    & 10   & 11   & 12   & 13   & 14   & 15   & 16   & 17   & m1   & m2   & D1   & D2   \\
\midrule
1  & 0   & 10.2 & 15.8 & 11.4 & 14.6 & 3.8  & 10.9 & 4.5  & 7.9  & 7.7  & 7.7  & 16.3 & 10.3 & 2.3  & 4.8  & 6.4  & 7.0    & 8.4  & 3.9  & 12.0   & 3.9  \\
2  & 9.1  & 0    & 20.8 & 17.7 & 23.2 & 13.0   & 15.1 & 7.8  & 12.0 & 13.6 & 11.9 & 20.6 & 1.1  & 11.3 & 11.4 & 15.4 & 8.7  & 12.6 & 12.1 & 20.6 & 11.1 \\
3  & 15.7  & 21.0    & 0    & 3.9  & 17.6 & 12.5 & 5.1  & 13.3 & 8.1  & 9.2  & 8.8  & 1.6  & 21.1 & 13.3 & 7.9  & 10.3 & 16.2 & 6.9 & 13.4 & 14.0 & 18.2 \\
4   & 11.3 & 17.2  & 3.9  & 0    & 13.2 & 8.1  & 8.3  & 10.6 & 8.3  & 6.5  & 8.1  & 4.8  & 17.3 & 8.9  & 6.1  & 5.9  & 12.9 & 4.6   & 9.0  & 9.6  & 13.8 \\
5  & 14.8  & 23.0 & 17.6 & 13.2 & 0    & 11.4 & 23   & 19.4 & 22.8 & 17.9 & 22.6 & 18.2 & 23.2 & 12.2 & 15.8 & 9.2  & 19.8 & 16.8 & 12.3 & 6.3  & 18.1 \\
6  & 4.5  & 12.7 & 12.6 & 8.3  & 11.5 & 0    & 14.8 & 8.5  & 11.8 & 8.3  & 11.7 & 13.3 & 12.9 & 2.6  & 5.8  & 3.0  & 8.9  & 8.1  & 0.8  & 8.0  & 7.9  \\
7  & 10.0 & 15.2   & 5.2  & 8.3  & 23.1 & 13.9 & 0    & 8.7  & 3.5  & 5.2  & 4.1  & 5.5  & 15.5 & 12.3 & 7.9  & 15.2 & 11.5 & 5.8 & 12.8 & 18.9 & 12.5 \\
8  & 5.1   & 7.9 & 10.7 & 10.4 & 19.2 & 9.0  & 8.2  & 0    & 5.2  & 4.1  & 4.8  & 13.7 & 8.0  & 7.4  & 3.5  & 10.5 & 4.1  & 4.8  & 6.4  & 16.6 & 7.6  \\
9  & 6.7  & 12.0 & 8.3  & 8.6  & 20.8 & 10.6 & 3.5  & 4.8  & 0    & 2.6  & 0.9  & 8.6  & 11.9 & 8.9  & 4.9  & 13.0 & 8.2  & 3.6  & 9.8  & 18.2 & 9.2  \\
10 & 6.1  & 11.3 & 9.2  & 7.2  & 19.1 & 10.0 & 5.2  & 4.4  & 2.3  & 0    & 2.1  & 9.5  & 11.4 & 7.8  & 2.9  & 11.5 & 7.6  & 1.5  & 7.7  & 15.4 & 8.6  \\
11 & 6.8  & 11.8 & 8.3  & 8.6  & 20.9 & 10.7 & 3.5  & 5.4  & 0.3  & 2.7  & 0    & 8.6  & 12.0 & 9.1  & 4.9  & 13.2 & 8.3  & 3.6  & 9.6  & 18.3 & 9.4  \\
12 & 15.0 & 20.3 & 1.7  & 4.8  & 18.5 & 13.4 & 5.5  & 13.6 & 8.6  & 9.6  & 9.2  & 0    & 20.4 & 14.5 & 8.8  & 11.2 & 16.5 & 7.8 & 14.3 & 14.8 & 17.5 \\
13 & 9.7  & 1.5  & 21.6 & 18.5 & 23.9 & 13.8 & 15.8 & 8.6  & 12.8 & 13.2 & 12.6 & 21.3 & 0    & 12.1 & 12.2 & 15.1 & 9.5  & 13.4 & 12.3 & 21.4 & 11.7 \\
14 & 3.4  & 11.6 & 13.2 & 8.8  & 12.0 & 2.6  & 12.3 & 5.9  & 9.3  & 7.1  & 9.1  & 13.8 & 11.8 & 0    & 4.6  & 3.1  & 8.4  & 6.9  & 2.5  & 9.2  & 6.8  \\
15 & 5.0  & 11.3 & 8.4  & 5.6  & 15.9 & 6.5  & 9.0  & 4.6  & 4.3  & 2.9  & 4.0  & 9.0  & 11.4 & 4.7  & 0    & 7.6  & 7.5  & 2.7  & 5.2  & 12.3 & 8.0  \\
16 & 7.2  & 15.4 & 10.3 & 5.9  & 9.1  & 2.9  & 14.7 & 8.4  & 11.7 & 11.5 & 11.5 & 10.9 & 14.2 & 3.3  & 7.1  & 0    & 10.8 & 9.4  & 3.8  & 6.3  & 10.5 \\
17 & 4.4  & 8.7  & 16.2 & 13.4 & 18.5 & 8.3  & 11.0 & 3.5  & 8.2  & 7.9  & 7.8  & 16.5 & 8.1  & 6.6  & 6.8  & 10.3 & 0    & 8.5  & 8.6  & 15.9 & 6.4  \\
m1 & 6.6  & 11.9 & 7.6 & 4.9  & 17.0 & 10.5  & 5.9 & 4.6  & 3.3 & 1.6  & 3.3 & 8.9 & 12.0 & 5.8    & 2.0    & 9.7  & 8.1  & 0    & -    & 13.3  & 9.1  \\
m2 & 3.9  & 12.4 & 14.4 & 8.9  & 12.2 & 0.9  & 13.4 & 8.5  & 10.1 & 9.4  & 10.2 & 14.5 & 13.2 & 3.0  & 6.0  & 3.7  & 7.7  & -    & 0    & 8.7  & 7.5  \\
D1 & 11.7 & 19.9 & 13.8 & 9.4  & 6.7  & 7.9  & 19.3 & 14.2 & 17.6 & 14.1 & 17.4 & 14.4 & 20.1 & 9.1  & 12.0 & 6.2  & 16.7 & 13.0  & 8.7  & 0    & -    \\
D2 & 3.6  & 10.9 & 17.8 & 13.8 & 18.0 & 7.7  & 12.0 & 5.7  & 9.0  & 10.6 & 8.9  & 17.5 & 11.1 & 6.1  & 7.2  & 9.8  & 7.1  & 9.6  & 7.2  & -    & 0   \\
\bottomrule
\end{tabular*}
\end{sidewaystable}

\bibliographystyle{cas-model2-names}

\bibliography{reference}

\end{document}